\documentclass[journal,twoside]{IEEEtran}
\usepackage{}

\usepackage{amsfonts}
\usepackage{amssymb}
\usepackage{amsmath}
\usepackage{amsmath,amssymb}
\usepackage{amsmath}
\usepackage{graphicx}
\usepackage{amsmath,amssymb,amsthm}
\usepackage{leftidx}
\usepackage{extarrows}
\usepackage{caption}
\usepackage{enumerate}
\usepackage{graphics}
\usepackage{subfigure}
\usepackage{float}
\usepackage{pict2e}
\usepackage{tikz}
\usepackage{color}
\usepackage{bm}

\usepackage[ruled,linesnumbered]{algorithm2e}  
\usepackage{amsmath}

\ifCLASSINFOpdf
\else
\fi
\begin{document}
%
\title{\textbf{Event-Triggered Optimal Attitude Consensus \\of Multiple Rigid Body Networks with \\Unknown Dynamics}}
%
%
%

\author{Xin Jin, Shuai Mao, Ljupco Kocarev, \IEEEmembership{Fellow,~IEEE}, Chen Liang, Saiwei Wang, \\and Yang Tang, \IEEEmembership{Senior Member,~IEEE}
\thanks{This work was supported in part by the National Key Research and Development Program of China under Grant 2018YFC0809302, the National Natural Science Foundation of China under Grant 61751305, in part by the Program of Shanghai Academic Research Leader under Grant No. 20XD1401300, in part by the Programme of Introducing Talents of Discipline to Universities (the 111 Project) under Grant B17017. \emph{(Corresponding
		author: Yang Tang.)}}

\thanks{Xin Jin, Shuai Mao, Saiwei Wang, Chen Liang and Yang Tang are with the Key Laboratory of Smart Manufacturing in Energy Chemical Process,
Ministry of Education, East China University of Science and Technology, Shanghai 200237, China (e-mail: jx\_9810@163.com; mshecust@163.com; 842376017@qq.com; tangtany@gmail.com; chenliang@ecust.edu.cn).}

\thanks{Ljupco Kocarev is with the Macedonian Academy of Sciences and Arts,
	1000 Skopje, Macedonia, also with the Faculty of Computer Science and
	Engineering, Univerzitet ``Sv. Kiril i Metodij," 1000 Skopje, Macedonia, and also with the BioCircuits Institute, University of California at San Diego, La Jolla, CA 92093 USA (e-mail: lkocarev@manu.edu.mk).}
 
}

\maketitle

\begin{abstract}
In this paper, \textcolor{black}{an event-triggered Reinforcement Learning (RL) method is  proposed for the optimal attitude consensus of multiple rigid body networks with unknown dynamics. }
Firstly, the consensus error is constructed through the attitude dynamics. According to the Bellman optimality principle, the  implicit form of the  optimal controller and the corresponding Hamilton-Jacobi-Bellman (HJB) equations are obtained. Because of the augmented system, the optimal controller can be obtained directly without relying on  the system dynamics.
Secondly, the  self-triggered mechanism is applied to reduce the computing and communication burden when updating the controller. 
In order to address the problem that the HJB equations are difficult to solve analytically, a RL method which only requires  measurement data at the event-triggered instants is proposed.
For each agent, only one neural network is designed to approximate the optimal value function.
\textcolor{black}{Each neural network is  updated only at the event-triggered instants. Meanwhile, the Uniformly
	Ultimately Bounded (UUB) of the closed-loop system is obtained, 
and the Zeno behavior is also avoided.} Finally, the simulation results on a multiple rigid body network demonstrate the validity of the proposed method. 
\end{abstract}


%
\IEEEpeerreviewmaketitle

\section{Introduction}\label{section1}
%
%
%
%


\IEEEPARstart{C}{onsensus} control, as a fundamental form of coordination problem in multi-agent systems, aims to design a control protocol for each agent to drive the states of all agents to be synchronized \cite{S.Ren2019}--\cite{D.Chen2020,J.Cao2019}.
Over recent decades, the attitude consensus problem of multiple rigid body networks has  received increasing attention  \cite{N.Xiong2017_sensors} because  it plays a significant role in the development of many fields, such as formation in three-dimensional space \cite{H.Du2019}, \cite{Z.Li2019}, cooperation of multi-manipulators \cite{X.Jin2019} and satellite networks \cite{H.Zhang2018}, etc. At present, some results have been proposed, which can be classified into two categories, the leaderless attitude consensus \cite{K.Zhang2013}--\cite{A.Abdessameud2009} and the leader-follower attitude consensus \cite{H.Cai2016}--\cite{M.Lu2020}. Note that none of them have considered the performance cost when achieving the attitude consensus.

In practical application scenarios, the performance cost is a factor that must be considered, which affects the efficiency of mission completion and the endurance of limited resources. 
\textcolor{black}{The optimal attitude consensus control, which not only makes the attitude of all rigid body systems tend to be synchronized, but also minimizes the performance cost.} In general, the optimal control problem can be transformed into solving the Hamilton-Jacobi-Bellman (HJB) equations. Nevertheless, it is very difficult to find the analytic solutions to the HJB equations. {\color{black}With the popularity of reinforcement learning technology \cite{B.Yi2019,B.Lin2019,N.Xiong2018} and the rapid development of the computing capacity of processors, some reinforcement learning based researches on solving the optimal consensus problem have emerged.} As far as we know, the vast majority of the research objects are linear systems \cite{K.G.Vamvoudakis2012}--\cite{J.Qin2019} or first-order nonlinear systems \cite{M.Abu-Khalaf2005}. In the above results, the knowledge of system dynamics is required in  \cite{K.G.Vamvoudakis2012}, \cite{M.Abu-Khalaf2005}, and the researches in \cite{J.Li2017} and  \cite{J.Qin2019} circumvent the dependence on system dynamics. However, the implementation of  algorithms in \cite{J.Li2017} and \cite{J.Qin2019} requires the acquisition of measurement data in advance and a lot of tedious integration operations are involved, which obviously increases the computational burden of the system \cite{N.Xiong2012}. At present, there are relatively few results concentrating on the reinforcement learning based method to realize the optimal attitude consensus for multiple rigid body networks. In \cite{H.Zhang2019}, a model-free algorithm is proposed to deal with the optimal consensus for multiple rigid body networks, in which the model of each rigid body is expressed in the form of Euler-Lagrange equation. However, an extra neural network-based observer is designed to estimate the system dynamics, which imposes additional computational burden. Motivated by these factors, we aim to design a method that only needs real-time measurement data to achieve the optimal attitude consensus of multiple rigid body networks with unknown dynamics.

Updating the controller and the neural networks at each sampling instant based on the reinforcement learning method will take up a lot of computing and communication resources, especially when the system scale is huge. Therefore,  it is particularly necessary to integrate the event-triggered mechanism into the reinforcement learning method to reduce the consumption of resources. {\color{black}In recent years, the event-triggered control scheme is widely studied to save the control cost and energy resources \cite{B.Li2018_TCYB,S.Zhu2020}. In \cite{L.Dong2017}--\cite{Q.Zhang2018}, the event-triggered mechanism is introduced to solve the optimal control of an individual system. The optimal consensus of multi-agent systems is considered in \cite{W.Zhao2019}--\cite{Z.Shi2020} by using the event-triggered reinforcement learning method. However, the event-triggered conditions in all of the above event-triggered reinforcement learning methods \cite{L.Dong2017}--\cite{Z.Shi2020} \textcolor{black}{include the continuous information.}} Therefore, all agents need to obtain the state information of themselves and their neighbors in real-time to determine whether the event-triggered condition is satisfied, \textcolor{black}{which  inevitably  increases the communication resources.} 
\textcolor{black}{
Inspired by \cite{A.Girard2015}--\cite{X.Jin2020}, we aim  to design an event-triggered reinforcement learning method under the self-triggered mechanism,  thereby greatly reducing the consumption of  computing and communication resources. 
Compared with the common linear system and first-order nonlinear system \cite{K.G.Vamvoudakis2012}--\cite{M.Abu-Khalaf2005}, it is challenging to combine the self-triggered mechanism with the reinforcement learning method to solve the optimal attitude consensus problem of multiple rigid body networks, since the dynamic model of a rigid body is a second order system with state coupled characteristics and the underlying attitude configuration space is non-Euclidean.
}

In this paper, we deal with the optimal attitude consensus problem for multiple rigid body networks with unknown system dynamics. The dynamic event-triggered  mechanism is first introduced, which can significantly reduce the consumption of computing resources caused by updating the controller. Based on the discussion of the dynamic event-triggered condition, a sufficient self-triggered condition is proposed. Under the self-triggered mechanism, the continuous communication between rigid bodies can be avoided.
Moreover, a  reinforcement learning method is used to obtain the optimal policy. In detail, a rigid body only needs a neural network to approximate the optimal value function because of the existence of the augmented system \cite{Y.Liu2019}. Each neural network is updated only when the self-triggered condition is violated. The main contributions  are as follows:

{\color{black}
1) By using only the measurement data at the event-triggered instants, we achieve the optimal attitude consensus of multiple rigid body networks with unknown system dynamics. No additional actor neural network \cite{J.Li2017}, \cite{J.Qin2019}  or additional neural network-based observer \cite{H.Zhang2019} is used in this paper, which obviously reduces the complexity of the algorithm implementation.

2) Compared with the results in \cite{H.Zhang2019} and \cite{S.Wang2020},  both a dynamic event-triggered condition and a self-triggered condition are  integrated into the proposed reinforcement learning based method to solve the optimal attitude consensus of multiple rigid body networks. Under the self-triggered  mechanism, the neural networks are updated only at the event-triggered instants. Meanwhile, the continuous communication is also avoided. Therefore, the consumption of computing and communication resources would be greatly reduced.}

The following layout of this paper is as follows. In Section \ref{section2}, we introduce the notations used in this paper and the basics of graph theory. The model-free optimal attitude consensus problem is described in Section \ref{section3}. Meanwhile, the event-triggered mechanism is also introduced. We design an event-triggered reinforcement learning method in Section \ref{section4}. The feasibility of this method is  verified through a simulation in Section \ref{section5}. Section \ref{section6} gives the conclusion. 

\section{Preliminaries}\label{section2}

\subsection{Notations}\label{section2.1}

Throughout this paper, $\mathbb{R}$  represents the set of all real numbers, $\mathbb{R}_{>0}$ represents the set of all positive real numbers,  $\mathbb{N}$ represents the set of all non-negative integers, and $\mathbb{N}_{>0}$ is the set of all positive integers, i.e., $\mathbb{R}=(-\infty, +\infty)$, $\mathbb{R}_{>0}=(0, +\infty)$, $\mathbb{N}=\{0,1,2,...\}$ and $\mathbb{N}_{>0}=\{1,2,...\}$.
$x \in \mathbb{R}^n$ indicates an $n$-dimensional vector, $I_n$ indicates an $n$-dimensional identity matrix, $A \in \mathbb{R}^{n\times m}$ indicates an $n\times m$ dimensional matrix.
For a vector $x$, its Euclidean norm is defined as $\lVert x\rVert=\sqrt{x^{T}x}$.
For a square matrix $B=[b_{ij}]\in \mathbb{R}^{n\times n}$, its trace is defined as $\textrm{tr}(B)=\sum_{i=1}^n b_{ii}$, its Frobenius norm is defined as $\lVert B\rVert =\sqrt{\sum_{i=1}^n \sum_{j=1}^n \vert b_{ij}\vert^2}$.  Define ${\lambda}_\textrm{min}(B)$ and ${\lambda}_\textrm{max}(B)$ as the minimum eigenvalue and maximum eigenvalue, respectively. $B>0$ $(B \ge 0)$ indicates that $B$ is positive (semi-positive) definite.  

For any two vectors $\xi=[x_1,y_1,z_1]^{\top} \in \mathbb{R}^3$ and $\zeta=[x_2,y_2,z_2]^{\top}\in \mathbb{R}^3$,  their cross product is expressed as follows:
\begin{align}\nonumber
\textcolor{black}{\xi \times \zeta}=
\begin{bmatrix}
y_1z_2-y_2z_1\\
x_2z_1-x_1z_2\\
x_1y_2-x_2y_1
\end{bmatrix} \in \mathbb{R}^3.
\end{align}

\subsection{Graph Theory}\label{section2.2}
Let $\mathcal{G}=(\mathcal{V},\mathcal{E})$ represent the directed communication graph among $N \in \mathbb{N}_{>0}$ rigid bodies, where $\mathcal{V}=\{1,2,...,N\}$ indicates the set of all rigid bodies and $\mathcal{E} \subseteq \mathcal{V} \times \mathcal{V}$ indicates  the communication relationship between any two rigid bodies. \textcolor{black}{For the rigid body $i$}, we use $\mathcal{N}_i=\{j\in\mathcal{V}|(j,i)\in\mathcal{E}\}$  to represent the set of its neighbors. For any two rigid bodies, if there is always a direct path between them, we call this communication graph strongly connected. \textcolor{black}{In this paper, we suppose that all directed communication graphs are strongly connected.}

In order to express the communication relationship between all rigid bodies more clearly, the weighted adjacency matrix $A=[a_{ij}] \in \mathbb{R}^{N\times N}$ is introduced.  If {\color{black} the rigid body $i$ can receive the data transmitted by the rigid body $j$, $a_{ij}>0$ and $a_{ij}=0$ otherwise.} 
The in-degree matrix of the directed communication graph can be expressed as $D=\text{diag}\{d_1,d_2,...,d_N\}\in \mathbb{R}^{N\times N}$, where $d_i=\sum_{j \in \mathcal{N}_i}a_{ij}$.
Let $\mathcal{L}=D-A=[l_{ij}]$ represent the Laplacian matrix, where $l_{ii}=d_i$, and $l_{ij}=- a_{ij}$ when $i \ne j$.

\section{Model-Free Event-Triggered Optimal Attitude Consensus}\label{section3}
\subsection{Model-Free Optimal Attitude Consensus}\label{section3.1}
{\color{black}We consider a multiple rigid body network with $N$ nodes, where the attitude of each node can be expressed by Modified Rodriguez Parameters (MRPs) \cite{H.Schaub2003}. }
{\color{black} For {the rigid body $i$}, the attitude is represent by $\sigma_i=[\sigma_i^1,\sigma_i^2,\sigma_i^3]^{\top}=\Psi_i\textrm{tan}\frac{\Phi_i}{4} \in \mathbb{R}^3$, where $\Psi_i\in \mathbb{R}^{3}$ indicates the Euler axis, and $\Phi_i\in [0, \pi)$ denotes the angle respective to the Euler axis.}

Then, the attitude dynamics of each rigid body is given in the following form:
\begin{subequations}\label{attitude dynamics}
	\begin{align}
	\dot{\sigma_i}=G(\sigma_i)\omega_i,&\\
	J_i\dot{\omega}_i=-\omega_i \times (J_i\omega_i)+\tau_i,& \quad i=1, 2, ...,N,
	\end{align}	
\end{subequations}
where $\omega_i \in \mathbb{R}^3$, $J_i \in \mathbb{R}^{3 \times 3}$  and $\tau_i \in \mathbb{R}^3$ indicate the angular velocity vector, the inertia matrix and the control input torque, respectively. The matrix $G(\sigma_i)=\frac{1}{2}(\sigma_i^{\times}+\sigma_i\sigma_i^{\top}+\frac{1-\sigma_i^{\top}\sigma_i}{2}I_3) \in \mathbb{R}^{3 \times 3}$, where 
\begin{align}
\nonumber
\sigma_i^{\times}=
\begin{bmatrix}
0&-\sigma_i^{3}&\sigma_i^{2}\\
\sigma_i^{3}&0&-\sigma_i^{1}\\
-\sigma_i^{2}&\sigma_i^{1}&0
\end{bmatrix}.
\end{align}

\textit{Definition 1:} \textcolor{black}{Given that the communication topology of a multiple rigid body network (\ref{attitude dynamics}) is strongly connected, the attitude consensus is said to be achieved  when the following conditions hold: }
\begin{subequations}\label{attitude consensus}
	\begin{align}
	&\lim_{t \to \infty} \big\lVert \sigma_i(t)-\sigma_j(t)\big\rVert=0, \\
	&\lim_{t \to \infty} \big\lVert \omega_i(t)-\omega_j(t) \big\rVert=0, \ \forall i,j \in \mathcal{V}.	
	\end{align}	
\end{subequations}

Considering the communication topology among these $N$ rigid bodies,  we can define the following form of consensus error \textcolor{black}{for the rigid body $i$}:
\begin{align} \label{consensus error}
\delta_i=\sum_{j\in \mathcal{N}_i}a_{ij}(\omega_i-\omega_j)+\alpha_i\sum_{j\in \mathcal{N}_i}a_{ij}(\sigma_i-\sigma_j),
\end{align}
where $\alpha_i \in \mathbb{R}_{>0}$.
When  $\lim_{t \to \infty}\delta_i=0$, $i=1,2,...,N$, we can easily obtain that  $\sigma_1=\sigma_2=...=\sigma_N$ and $\omega_1=\omega_2=...=\omega_N$ with $t\to \infty$. That is to say, the attitude consensus is achieved.

The dynamics of $\delta_i$ can be obtained by taking the derivative of  Eq. (\ref{consensus error}), which is described as follows:
\begin{align}\label{dot consensus error}
\nonumber\dot{\delta}_i&=\sum_{j\in \mathcal{N}_i}a_{ij}(\dot{\omega}_i-\dot{\omega}_j)+\alpha_i\sum_{j\in \mathcal{N}_i}a_{ij}(\dot{\sigma}_i-\dot{\sigma}_j)\\
& =\varGamma_i(\delta_i)+l_{ii}J_i^{-1}\tau_i-\sum_{j\in \mathcal{N}_i}a_{ij}J_j^{-1}\tau_j,
\end{align}
where $ \varGamma_i(\delta_i)=\alpha_i\sum_{j\in \mathcal{N}_i}a_{ij}(\dot{\sigma}_i-\dot{\sigma}_j)+\sum_{j\in \mathcal{N}_i}a_{ij}\Big[-J_i^{-1}\Big(\big(G^{-1}(\sigma_i)\dot{\sigma_i}\big) \times \big(J_iG^{-1}(\sigma_i)\dot{\sigma_i}\big)\Big)
+J_j^{-1}\Big(\big(G^{-1}(\sigma_j)\dot{\sigma_j}\big) \times \big(J_jG^{-1}(\sigma_j)\dot{\sigma_j}\big)\Big)\Big]$.

In order to overcome the dependence on model information, a compensator is introduced, which can be expressed by the following affine differential equation:
\begin{align}\label{compensator}
\dot{\tau}_i=f(\tau_i)+l_{ii}g(\tau_i)u_i-\sum_{j \in \mathcal{N}_i}a_{ij}g(\tau_j)u_j,
\end{align}
where $f(\cdot): \mathbb{R}^3 \to \mathbb{R}^3$, $g(\cdot): \mathbb{R}^3 \to \mathbb{R}^{3 \times 3}$ are two functions to be designed later, and $u_i \in \mathbb{R}^3$ is the control input of the compensator. We need to choose appropriate functions $f(\cdot)$ and $g(\cdot)$ to ensure that  the compensator is controllable. In this paper, a feasible pair of $f(\cdot)$ and $g(\cdot)$ is given as follows:
\begin{subequations}
	\begin{align}
		&f(\tau_i)=-2\tau_i,\\
		&g(\tau_i)=\text{diag}\{\cos^2(\tau_i^1),\cos^2(\tau_i^2),\cos^2(\tau_i^3)\}.
	\end{align}	
\end{subequations}

By combining the consensus error $\delta_i$ and the control input torque $\tau_i$, we define an augmented consensus error, which is expressed as $e_i=[\delta_i^{\top},\tau_i^{\top}]^{\top} \in \mathbb{R}^6$. According to Eq. (\ref{dot consensus error}) and Eq. (\ref{compensator}), we can use the following augmented system to describe the dynamics of the augmented consensus error:
\begin{align}\label{augment system}
\dot{e}_i=X_i(e_i)+l_{ii}Y_i(e_i)u_i-\sum_{j\in \mathcal{N}_i}a_{ij}Y_j(e_j)u_j,
\end{align}
where $X_i(e_i)$ and $Y_i(e_i)$ are represented as follows:
\begin{align}
\nonumber X_i(e_i)=\begin{bmatrix}
\varGamma_i(\delta_i)+l_{ii}J_i^{-1}\tau_i-\sum\limits_{j\in \mathcal{N}_i}a_{ij}J_j^{-1}\tau_j\\
f(\tau_i)
\end{bmatrix} \in \mathbb{R}^6,
\end{align}

\begin{align}
\nonumber Y_i(e_i)=\begin{bmatrix}
0 \\
g(\tau_i)
\end{bmatrix} \in \mathbb{R}^{6\times 3}.
\end{align}

\textit{Assumption 1:} The matrix $X_i(e_i)$ is bounded, i.e., $\forall i \in \mathcal{V}$,  $\lVert X_i(e_i) \rVert \leq X_M \lVert e_i \rVert$ is satisfied, where $X_M \in \mathbb{R}_{>0}$.

In order to measure the performance cost of implementing the attitude consensus, a performance function is defined in the following form:
\begin{align}\label{performance function}
F_i(e_i(0),u_i,u_{-i})=\int_{0}^\infty\big(e_i^{\top}Q_ie_i+u_i^{\top}R_iu_i\big) dt,
\end{align}
where $u_{-i}=\{u_j|j\in \mathcal{N}_i\}$ indicates the set of control inputs for the neighbors of \textcolor{black}{the rigid body $i$}, $Q_i\in \mathbb{R}^{6 \times 6}$, $Q_i\ge 0$, $R_i\in \mathbb{R}^{3 \times 3}$ and $R_i>0$. \textcolor{black}{According to Eq. (\ref{augment system}), we can conclude that $e_i$ is driven by $u_i$ and $u_{-i}$. Therefore, the left side of  Eq. (\ref{performance function}) also contains $u_{-i}$.}

According to (\ref{performance function}), the value function can be defined as
\begin{align}\label{value function}
V_i(e_i(t))=\int_{t}^\infty\big(e_i(v)^{\top}Q_ie_i(v)+u_i(v)^{\top}R_iu_i(v)\big) dv.
\end{align}

By taking the derivative of Eq. (\ref{value function}), we can obtain the  Hamiltonian function in the following form:
\begin{align}\label{Hamiltonian}
\nonumber   &H_i(e_i,\nabla V_i, u_i, u_{-i}) \\ 
\nonumber   &=e_i^{\top}Q_ie_i+u_i^{\top}R_iu_i+\nabla V_i^{\top}\Big(X_i+l_{ii}Y_iu_i-\sum_{j \in \mathcal{N}_i}a_{ij}Y_ju_j\Big)\\ 
&=0,
\end{align}
where $\nabla V_i=\frac{\partial V_i}{\partial e_i}.$

The implicit solution of the model-free optimal controller $u_i^*$ can be derived from $\frac{\partial H_i}{\partial u_i}=0 $, which is represented as
\begin{align}\label{optimal controller}
u_i^*=-\frac{1}{2}l_{ii}R_i^{-1} Y_i^{\top} \nabla V_i^*,
\end{align}
where $V_i^*$ indicates the optimal value function and $\nabla V_i^*=\frac{\partial V_i^*}{\partial e_i}$. 

Combining (\ref{Hamiltonian}) and (\ref{optimal controller}), we can derive the HJB equation for \textcolor{black}{the rigid body $i$}  as follows:
\begin{align}\label{HJB}
\nonumber   &H_i(e_i,\nabla V_i^*, u_i^*, u_{-i}^*) \\ 
\nonumber   &=e_i^{\top}Q_ie_i+(u_i^*)^{\top}R_iu_i^*\\
\nonumber &\quad+ (\nabla V_i^*)^{\top} \Big(X_i+l_{ii}Y_iu_i^*-\sum_{j \in \mathcal{N}_i}a_{ij}Y_ju_j^*\Big)\\ 
\nonumber   &=e_i^{\top}Q_ie_i-\frac{1}{4}l_{ii}^2(\nabla V_i^*)^{\top}Y_iR_i^{-1}Y_i^{\top}\nabla V_i^*\\
\nonumber &\quad+(\nabla V_i^*)^{\top}\Big(X_i+\frac{1}{2}\sum_{j \in \mathcal{N}_i}a_{ij}l_{jj}Y_jR_j^{-1}Y_j^{\top}\nabla V_j^*\Big)\\
&=0.
\end{align}

From Eq. (\ref{optimal controller}), we can observe that the optimal controller contains $Y_i$ and $V_i^*$.
According to the definition of the augmented system (\ref{augment system}), we can obtain $Y_i$ which overcomes the dependence on system dynamics. Therefore, we only need to obtain the optimal value function $V_i^*$ from Eq. (\ref{HJB}) to achieve the optimal controller.

\subsection{Dynamic Event-Triggered Mechanism}
For the purpose of  reducing the computational burden when updating the controller, the event-triggered mechanism is introduced. \textcolor{black}{Under the event-triggered mechanism, we only update the controller at a series of discrete instants $\{t_i^h\}_{h \in \mathbb{N}}$, where $t_i^h<t_i^{h+1}$ holds with $\forall h \in \mathbb{N}$ and  the initial instant is set as $t_i^0=0$. } 

To define the difference between the augmented consensus error at the last event-triggered instant and in real-time  as the  measurement error $E_i(t) \in \mathbb{R}^6$, which is expressed as
\begin{align}\label{event-triggered consensus error}
E_i(t)= e_i(t_i^h)-e_i(t),\  \forall t \in [t_i^h,t_i^{h+1}).  
\end{align}

For \textcolor{black}{the rigid body $i$}, the controller $u_i$ is only updated at the event-triggered instants $t_i^h$, and remains unchanged until a new event is triggered. \textcolor{black}{During the event-triggered intervals $[t_i^h,t_i^{h+1})$,  the neighbor $j$ of the rigid body $i$ might update its controller at some instants $t_j^{h'}$, which performs  as a piecewise constant. Letting $\mu_j$ indicate the event-triggered times of  the neighbor $j$, we have $t_j^{h'} \in \{t_j^0,t_j^1,...,t_j^{\mu_j}\}$ and $t_j^0 = t_i^h$.}

Therefore, the $e_i$-dynamics (\ref{dot consensus error}) during $[t_i^h,t_i^{h+1})$ is represented as follows:
\begin{align}\label{event dot consensus error}
\dot{e}_i=X_i+l_{ii}Y_i\hat{u}_i-\sum_{j \in \mathcal{N}_i}a_{ij}Y_j\hat{u}_j,
\end{align}
where  $\hat{u}_i=u_i(t_i^h)$ and $\hat{u}_j=u_j(t_j^{h'})$ are the control input vectors at the event-triggered instants.

\textit{Definition 2 (Event-Triggered Admissible Control):} If the following conditions are  met: $u_i$ is
piecewise continuous, $u_i(0)=0$, $e_i$-dynamics (\ref{event dot consensus error})  is stable, and the performance function (\ref{performance function}) is finite, then $u_i$ is called event-triggered admissible control.

Depending on the form of optimal controller (\ref{optimal controller}),  we can obtain the event-triggered optimal controller as follows:
\begin{align}\label{event optimal controller}
\nonumber\hat{u}_i^*&=u_i^*(t_i^h)\\
&=-\frac{1}{2} l_{ii} R_i^{-1} Y_i^{\top} \nabla \hat{V}_i^*,\ \forall t \in [t_i^h,t_i^{h+1}),
\end{align}
where $\nabla \hat{V}_i^*=\frac{\partial V_i^*}{\partial e_i}(t_i^h)$.

By combining (\ref{HJB}) and (\ref{event optimal controller}), the event-triggered HJB equation for \textcolor{black}{the rigid body $i$} is given as follows:
\begin{align}\label{event HJB}
\nonumber   &H_i(e_i,\nabla \hat{V}_i^*, \hat{u}_i^*, \hat{u}_{-i}^*) \\ 
\nonumber &=e_i^{\top}Q_ie_i-\frac{1}{4}(\nabla \hat{V}_i^*)^{\top}Y_iR_i^{-1}Y_i^{\top}\nabla \hat{V}_i^*\\
\nonumber &\quad+(\nabla \hat{V}_i^*)^{\top}\Big(X_i+\frac{1}{2}\sum_{j \in \mathcal{N}_i}a_{ij}l_{jj}Y_jR_j^{-1}Y_j^{\top}\nabla \hat{V}_j^*\Big)\\
&=0.
\end{align}

The following assumption is proposed for proving the stability of system (\ref{event dot consensus error}).

\textcolor{black}{
\textit{Assumption 2 \cite{M.Lemmon2010}:} The controller $u_i$ is Lipschitz continuous during the time interval $[t_i^h, t_i^{h+1})$, and there exists a constant $P$ that satisfies the following inequality:
\begin{align}\label{Lipschitz}
\big\lVert u_i\big(e_i(t_i^h)\big)-u_i\big(e_i(t)\big)\big\rVert \leq P\lVert E_i(t)\rVert,
\end{align}
where $P$ indicates the Lipschitz constant. In the actual engineering applications, the selection of $P$ should not be smaller than the maximum value of $\lVert \partial u_i/\partial e_i^{\top} \rVert$.
}

For the event-triggered control, it is necessary to ensure that Zeno-behavior does not occur in the proposed  event-triggered mechanism. Therefore, we introduce a dynamic  event-triggered mechanism  inspired by \cite{A.Girard2015} to exclude the Zeno-behavior implicitly. Firstly, a dynamic  variable $y_i(t)$ is defined as follows:
\begin{align}\label{y}
\nonumber\dot{y}_i(t)=&-\gamma_i y_i(t)+ \kappa_i\Big( \varpi_i \lambda_\textrm{min}(Q_i) \big\lVert e_i(t) \big\rVert^2\\
&-\lambda_\textrm{max}(R_i)P^2 \big\lVert E_i(t) \big\rVert^2 \Big),
\end{align}
where $y_i(0) \in \mathbb{R}_{\ge 0}$, $\gamma_i \in \mathbb{R}_{> 0}$, $\kappa_i \in [0,\frac{1}{2}]$ and $\varpi_i \in [0,1]$. Then, we can obtain the event-triggered instants through the following dynamic event-triggered condition:
\begin{align}\label{dynamic ETC}
\nonumber t_i^0=0,\\
\nonumber t_i^{h+1}=&\max\limits_{r \ge t_i^h}\biggl\{r\in\mathbb{R}:y_i(t)+\theta_i\Big( \varpi_i \lambda_\textrm{min}(Q_i) \big\lVert e_i(t) \big\rVert^2\\
&-\lambda_\textrm{max}(R_i)P^2 \big\lVert E_i(t) \big\rVert^2 \Big)\ge 0, \forall t \in [t_i^h,r]\biggr\},
\end{align}
where $\theta_i \in \mathbb{R}_{>0}$  will be determined later in the proof analysis.

\textit{Lemma 1:} Assuming that the event-triggered instants $t_i^h$ are determined by (\ref{dynamic ETC}), $y_i(t) \ge 0$  always holds for $\forall t \in [0,+\infty)$ if the  given initial value $y_i(0) \ge 0$.

\textit{Proof}: For $\forall t \in [0,+\infty)$, the  event-triggered condition (\ref{dynamic ETC}) guarantees the following  inequality:
\begin{align}\label{lamma1-1}
y_i(t)+\theta_i\Big( \varpi_i \lambda_\textrm{min}(Q_i) \big\lVert e_i(t) \big\rVert^2-\lambda_\textrm{max}(R_i)P^2 \big\lVert E_i(t) \big\rVert^2 \Big)\ge 0.
\end{align}

Since the selection of $\theta_i$ must satisfy $\theta_i>0$, the inequality (\ref{lamma1-1}) becomes
\begin{align}\label{lamma1-2}
\varpi_i \lambda_\textrm{min}(Q_i) \big\lVert e_i(t) \big\rVert^2-\lambda_\textrm{max}(R_i)P^2 \big\lVert E_i(t) \big\rVert^2 \ge -\frac{1}{\theta_i}y_i(t).
\end{align}

By combining  (\ref{y}) and (\ref{lamma1-2}), we can easily obtain that for $\forall t \in  [0,+\infty)$,
\begin{align}\label{lamma1-3}
\dot{y}_i(t) \ge -(\gamma_i+\frac{\kappa_i}{\theta_i})y_i(t).
\end{align}

According to the comparison lemma in \cite{H.K.Khalil2002}, we can  deduce
that
\begin{align}\label{lamma1-4}
y_i(t)\ge y_i(0) \textcolor{black}{\exp\{-(\gamma_i+\frac{\kappa_i}{\theta_i})t\}.}
\end{align}

Therefore, $y_i(t) \ge 0$ always holds for $\forall t \in [0,+\infty)$. The complete proof is given. $\hfill\blacksquare$

{\color{black}\textit{Theorem 2:} Consider a multiple rigid body network with $N$ nodes under a strongly connected communication topology. 
Supposed that Assumption 2 holds, the performance function and the event-triggered optimal controller are given by (\ref{performance function})  and (\ref{event optimal controller}), respectively.  
 If  the event-triggered instants   are determined by the dynamic event-triggered condition (\ref{dynamic ETC}), then the following two conclusions can be obtained:

1) The $e_i$-dynamics (\ref{event dot consensus error}) is asymptotically stable, i.e., the optimal attitude consensus is achieved.

2) The Zeno behavior is excluded, i.e., the interval between  $t_i^{h+1}$ and $t_i^h$, $\forall i \in \mathcal{V}$ has a positive lower bound.}

\textit{Proof}: 1) Firstly, we prove that the $e_i$-dynamics (\ref{event dot consensus error}) is asymptotically stable.
Choosing $\Pi_i(t)=V_i^*\big(e_i(t)\big)+y_i(t)$ as the Lyapunov function, which contains the optimal value function $V_i^*\big(e_i(t)\big)$ in (\ref{value function}) and the dynamic variable $y_i(t)$ is governed by (\ref{y}).

By taking the first-order derivative of $V_i^*(e_i)$ with respect to $t$ along with the  consensus error $e_i$,  we derive
\begin{align} \label{proof1-1}
\nonumber \dot{V}_i^*(e_i)&=(\nabla V_i^*)^{\top} \dot{e}_i\\
&=(\nabla V_i^*)^{\top}\Big(X_i+l_{ii}Y_i\hat{u}_i^*-\sum_{j\in \mathcal{N}_i}a_{ij}Y_j\hat{u}_j^*\Big).
\end{align}

During the event-triggered intervals $[t_i^h,t_i^{h+1})$, assuming that the neighbors of \textcolor{black}{the rigid body $i$} execute $\hat{u}_j^*=u_j^*(t_j^{h'})$. According to (\ref{optimal controller}) and (\ref{HJB}), it can be easily obtained that
\begin{align} \label{proof1-2}
\nonumber(\nabla V_i^*)^{\top}X_i=&-e_i^{\top}Q_ie_i-(u_i^*)^{\top}R_iu_i^*
-(\nabla V_i^*)^{\top}l_{ii}Y_iu_i^*\\
&+(\nabla V_i^*)^{\top}\sum_{j\in \mathcal{N}_i}a_{ij}Y_j\hat{u}_j^*
\end{align}
and
\begin{align}\label{proof1-3}
(\nabla& V_i^*)^{\top}l_{ii}Y_i=-2(u_i^*)^{\top}R_i.
\end{align}

Thus, Eq. (\ref{proof1-1}) can be redescribed as
\begin{align}\label{proof1-4}
\nonumber\dot{V}_i^*(e_i)=&-e_i^{\top}Q_ie_i-(u_i^*)^{\top}R_iu_i^*\\
\nonumber &+(\nabla V_i^*)^{\top}l_{ii}Y_i(\hat{u}_i^*-u_i^*)\\
\nonumber=&-e_i^{\top}Q_ie_i+(u_i^*)^{\top}R_iu_i^*-2(u_i^*)^{\top}R_i\hat{u}_i^*\\
\nonumber=&-e_i^{\top}Q_ie_i-(\hat{u}_i^*)^{\top}R_i\hat{u}_i^*\\
\nonumber&+(u_i^*-\hat{u}_i^*)^{\top}R_i(u_i^*-\hat{u}_i^*)\\
\leq&-\lambda_\textrm{min}(Q_i)\big\lVert e_i(t)\big\rVert^2+\lambda_\textrm{max}(R_i)P^2\big\lVert E_i(t)\big\rVert^2.
\end{align}

According to (\ref{y}) and (\ref{proof1-4}), we can obtain the first-order derivative of $\Pi_i(t)$ as follows:
\begin{align}\label{proof1-5} 
\nonumber \dot{\Pi}_i(t)=&\dot{V}_i^*(e_i)+\dot{y}_i(t)\\
\nonumber \leq& -\lambda_\textrm{min}(Q_i)\big\lVert e_i(t)\big\rVert^2+\lambda_\textrm{max}(R_i)P^2\big\lVert E_i(t) \big\rVert^2\\
\nonumber&-\gamma_i y_i(t)+ \kappa_i\Big( \varpi_i \lambda_\textrm{min}(Q_i) \big\lVert e_i(t) \big\rVert^2\\
\nonumber&-\lambda_\textrm{max}(R_i)P^2 \big\lVert E_i(t) \big\rVert^2 \Big)\\
\nonumber\leq&-(1- \varpi_i)\lambda_\textrm{min}(Q_i)\big\lVert e_i(t) \big\rVert^2-\gamma_i y_i(t)\\
\nonumber&+(\kappa_i-1)\Big( \varpi_i \lambda_\textrm{min}(Q_i) \big\lVert e_i(t) \big\rVert^2\\
&-\lambda_\textrm{max}(R_i)P^2 \big\lVert E_i(t) \big\rVert^2 \Big).
\end{align}

Substituting the dynamic event-triggered condition (\ref{dynamic ETC}) into (\ref{proof1-5}), $\dot{\Pi}_i(t)$ becomes
\begin{align}\label{proof1-6} 
\nonumber \dot{\Pi}_i(t) \leq& -(1- \varpi_i)\lambda_\textrm{min}(Q_i)\big\lVert e_i(t) \big\rVert^2-\gamma_i y_i(t)\\
\nonumber&+(\kappa_i-1)(-\frac{1}{\theta_i})y_i(t)\\
\leq&-(1- \varpi_i)\lambda_\textrm{min}(Q_i)\big\lVert e_i(t) \big\rVert^2-(\gamma_i+\frac{\kappa_i-1}{\theta_i})y_i(t).
\end{align}
Since  $\varpi_i \in [0,1]$, we have $\dot{\Pi}_i(t) \leq 0$ if $\theta_i \in [\frac{1-\kappa_i}{\gamma_i},+\infty)$.
Therefore, we can select appropriate $\gamma_i \in \mathbb{R}_{>0}$, $\kappa_i \in [0,\frac{1}{2}]$, $\varpi_i \in [0,1]$ and  $\theta_i \in [\frac{1-\kappa_i}{\gamma_i},+\infty)$ to ensure that the $e_i$-dynamics (\ref{event dot consensus error}) is asymptotically stable under the dynamic event-triggered condition (\ref{dynamic ETC}).

2) Then, we prove that the Zeno behavior is excluded.

According to (\ref{lamma1-4}), we can  deduce a sufficient condition of the dynamic event-trigger condition (\ref{dynamic ETC}) when $\varpi$ is selected as zero, which is expressed as follows:
{\color{black}\begin{align}\label{proof1-7} 
\nonumber t_i^{h+1}=&\max\limits_{r \ge t_i^h}\biggl\{r\in\mathbb{R}:\big\lVert E_i(t) \big\rVert \leq \sqrt{\frac{y_i(0)}{\theta_i\lambda_\textrm{max}(R_i)P^2}}\\
& \times \textcolor{black}{\exp\{-\frac{1}{2}{(\gamma_i+\frac{\kappa_i}{\theta_i})t}\}},\forall t \in [t_i^h,r]\biggr\}.
\end{align} }

According to the definition of $Y_i$,  we can conclude that $Y_i$ is bounded. That is to say, $Y_i \leq Y_M$ is  satisfied, where $Y_M \in \mathbb{R}_{>0}$.   With  Assumption 1 and the definition of measurement error $E_i(t)$, we can obtain that for $\forall t \in [t_i^h,t_i^{h+1})$,
\begin{align}\label{proof1-8}
\nonumber\big\lVert \dot{E}_i(t) \big\rVert&=\big\lVert \dot{e}_i(t_i^h)-\dot{e}_i(t) \big\rVert\\
\nonumber&=\big\lVert X_i\big(e_i(t)\big)+l_{ii}Y_iu_i(t_i^h)-\sum_{j\in \mathcal{N}_i}a_{ij}Y_ju_j(t_j^{h'})\big\rVert\\
\nonumber&\leq X_M\big\lVert e_i(t)\big\rVert+ \sum_{j\in \mathcal{N}_i}a_{ij}Y_M \big(\big\lVert u_i(t_i^h) \big \lVert+\big\lVert u_j(t_j^{h'}) \big\rVert\big)\\
\nonumber&\leq X_M\big\lVert E_i(t)\big\rVert+ X_M\big\lVert e_i(t_i^h)\big\rVert\\
&\quad+ \sum_{j\in \mathcal{N}_i}a_{ij}Y_M \big(\big\lVert u_i(t_i^h) \big \lVert+\big\lVert u_j(t_j^{h'}) \big\rVert\big).
\end{align} 

Since  $E_i(t_i^h)=0$ is satisfied at the event-triggered instants, we can derive the following inequality by using the comparison lemma \cite{H.K.Khalil2002}:
{\color{black}\begin{align}\label{proof1-9}
\nonumber \big\lVert E_i(t) \big\rVert &\leq \textcolor{black}{\exp\{X_M(t-t_i^h)\}}E_i(t_i^h)\\
\nonumber&\quad+\frac{1}{2}\int_{t_i^{h}}^t \textcolor{black}{\exp\{X_M(t-v)\}}\Big(X_M\big\lVert e_i(t_i^h)\big\rVert\\
\nonumber&\quad+ \sum_{j\in \mathcal{N}_i}a_{ij}Y_M \big(\big\lVert  u_i(t_i^h) \big \lVert+\big\lVert u_j(t_j^{h'}) \big\rVert\big)\Big) dv\\
\nonumber&\leq\frac{1}{2}\int_{t_i^{h}}^t \textcolor{black}{\exp\{X_M(t-v)\}}\Big(X_M\big\lVert e_i(t_i^h)\big\rVert\\
&\quad+ \sum_{j\in \mathcal{N}_i}a_{ij}Y_M \big(\big\lVert  u_i(t_i^h) \big \lVert+\big\lVert u_j(t_j^{h'}) \big\rVert\big)\Big) dv.
\end{align}}

Let $\tilde{t}_i^{h+1}$ indicate the next event-triggered instant determined by the sufficient condition (\ref{proof1-7}). According to (\ref{proof1-9}), the sufficient condition (\ref{proof1-7}) can be divided into two situations during the time interval $[t_i^h, \tilde{t}_i^{h+1})$, which contains: 1) there is no events occuring for all rigid bodies in $\mathcal{N}_i$, and 2) there exists at least one event  \textcolor{black}{for the rigid body $j\in \mathcal{N}_i$.}

\emph{Situation 1:} For all rigid bodies in $\mathcal{N}_i$, there is no event-triggered instants during $[t_i^h, \tilde{t}_i^{h+1})$.  Therefore, we can deduce the following inequality:
{\color{black}\begin{align}\label{proof1-10}
\nonumber&\sqrt{\frac{y_i(0)}{\theta_i\lambda_\textrm{max}(R_i)P^2}}\textcolor{black}{\exp\{-\frac{1}{2}{(\gamma_i+\frac{\kappa_i}{\theta_i})\tilde{t}_i^{h+1}}\}}\\
\nonumber &\leq 
\frac{X_M\big\lVert e_i(t_i^h)\big\rVert+ \sum\limits_{j\in \mathcal{N}_i}a_{ij}Y_M \big(\big\lVert u_i(t_i^h) \big \lVert+\big\lVert u_j(t_j^{h'}) \big\rVert\big)}{2X_M}\\
&\quad  \times \Big(\textcolor{black}{\exp\{X_M(\tilde{t}_i^{h+1}-t_i^h)\}}-1\Big).
\end{align}}

\emph{Situation 2:} \textcolor{black}{For the rigid body $j\in \mathcal{N}_i$}, there exists $\mu_j \in \mathbb{N}_{>0}$ event-triggered instants during  $[t_i^h, \tilde{t}_i^{h+1})$. By using $t_j^0, t_j^1, ..., t_j^{\mu_j}$  to indicate the event-triggered instants and $t_j^0=t_i^k$, we can deduce
{\color{black}\begin{align}\label{proof1-11}
\nonumber&\sqrt{\frac{y_i(0)}{\theta_i\lambda_\textrm{max}(R_i)P^2}}\textcolor{black}{\exp\{-\frac{1}{2}{(\gamma_i+\frac{\kappa_i}{\theta_i})\tilde{t}_i^{h+1}}\}}\\
\nonumber&\leq \frac{X_M\big\lVert e_i(t_i^h)\big\rVert+ \sum\limits_{j\in \mathcal{N}_i}a_{ij}Y_M \big\lVert  u_i(t_i^h) \big \lVert}{2X_M}\\
\nonumber&\quad\quad \times\Big(\textcolor{black}{\exp\{X_M(\tilde{t}_i^{h+1}-t_i^h)\}}-1\Big)\\
\nonumber&\quad +\frac{ \sum\limits_{j\in \mathcal{N}_i}a_{ij}Y_M \sum\limits_{s=0}^{\mu_j-1} \big\lVert  u_j(t_j^s) \big \lVert}{2X_M}\\
\nonumber&\quad\quad \times\Big(\textcolor{black}{\exp\{X_M(t_j^{s+1}-t_i^s)\}}-1\Big)\\
\nonumber &\quad +\frac{ \sum\limits_{j\in \mathcal{N}_i}a_{ij}Y_M  \big\lVert  u_j(t_j^{\mu_j}) \big \lVert}{2X_M}\\
&\quad\quad \times\Big(\textcolor{black}{\exp\{X_M(\tilde{t}_i^{h+1}-t_j^{\mu_j})\}}-1\Big).
\end{align}}

Combining (\ref{proof1-10}) and (\ref{proof1-11}), we can obtain the unified form of the two situations, which is expressed as follows: 
{\color{black}\begin{align}\label{proof1-12}
\nonumber&\sqrt{\frac{y_i(0)}{\theta_i\lambda_\textrm{max}(R_i)P^2}}\textcolor{black}{\exp\{-\frac{1}{2}{(\gamma_i+\frac{\kappa_i}{\theta_i})\tilde{t}_i^{h+1}}\}}\\
&\leq \frac{X_M\big\lVert e_i(t_i^h)\big\rVert+\Theta_i }{2X_M}\Big(\textcolor{black}{\exp\{X_M(\tilde{t}_i^{h+1}-t_i^h)\}}-1\Big),
\end{align}}
where $\Theta_i=\sum\limits_{j\in \mathcal{N}_i}a_{ij}Y_M \Big(\big\lVert  u_i(t_i^h) \big \lVert+\max\limits_{s=0,...,\mu_j}\big\{\big\lVert u_j(t_j^s) \big\rVert\big\}\Big)$.

Since $\tilde{t}_i^{h+1}$ is determined by the sufficient condition (\ref{proof1-7}), and let $t_i^{h+1}$ indicate the next event-triggered instant  determined by (\ref{dynamic ETC}), we can obtain the interval between two adjacent event-triggered instants:
{\color{black}\begin{align}\label{proof1-13}
\nonumber t_i^{h+1}-t_i^h &\ge \tilde{t}_i^{h+1}-t_i^h\\
\nonumber &\ge \frac{1}{X_M} \text{log}\Bigg(\frac{2X_M\sqrt{y_i(0)}}{\big(X_M\big\lVert e_i(t_i^h)\big\rVert+\Theta_i\big)\sqrt{\theta_i\lambda_\textrm{max}(R_i)P^2}}\\
&\quad\times \textcolor{black}{\exp\{-\frac{1}{2}{(\gamma_i+\frac{\kappa_i}{\theta_i})\tilde{t}_i^{h+1}}\}}+1\Bigg)>0.
\end{align}}

Therefore, the Zeno behavior can be excluded. The complete proof is given. $\hfill\blacksquare$ 

\subsection{Self-Triggered Mechanism}

Under the dynamic event-triggered mechanism, we have to obtain the continuous consensus error $e_i(t)$ and the continuous measurement error $E_i(t)$ to judge whether the dynamic event-triggered condition (\ref{dynamic ETC}) is violated. Therefore, it is necessary to continuously communicate with neighbors to obtain their absolute attitude information,  or to measure continuous relative attitude information with the help of sensors such as cameras.
In order to overcome this problem, a self-triggered condition is proposed in this subsection.

Letting $\kappa_i=0$, Eq. (\ref{y}) becomes
\begin{align}\label{self1}
\dot{y}_i(t)=-\lambda_iy_i(t).
\end{align}
Therefore, we can obtain that $y_i(t)=y_i(0)\textcolor{black}{\exp\{-\lambda_it\}}$, $\forall t \in [0,+\infty)$.

Then, letting $\varpi_i=0$, the dynamic event-triggered condition (\ref{dynamic ETC}) becomes
{\color{black}\begin{align}\label{self2}
\nonumber t_i^{h+1}=&\max\limits_{r \ge t_i^h}\biggl\{r\in\mathbb{R}: \big\lVert E_i(t) \big\rVert \leq \sqrt{\frac{y_i(0)}{\theta_i\lambda_\textrm{max}(R_i)P^2}}\\
&\times \textcolor{black}{\exp\{-\frac{1}{2}\lambda_it\}}, \forall t \in [t_i^h,r]\biggr\},
\end{align}}
where $\theta_i \in [\frac{1}{\gamma_i},+\infty)$.

According to (\ref{proof1-12}), the self-triggered measurement error is defined in the following form:
\begin{align}\label{self3}
\big\lVert \Delta_i(t) \big \rVert=\frac{X_M\big\lVert e_i(t_i^h)\big\rVert+\Theta_i }{2X_M}\Big(\textcolor{black}{\exp\{X_M(t-t_i^h)\}}-1\Big),
\end{align}
which is the upper bound of  $\big \lVert E_i(t)\big \rVert$. 

Thus, we can obtain a new sufficient condition of the dynamic event-triggered condition (\ref{dynamic ETC}) as follows:
{\color{black}\begin{align}\label{self4}
\nonumber t_i^{h+1}=&\max\limits_{r \ge t_i^h}\biggl\{r\in\mathbb{R}: \big\lVert \Delta_i(t) \big\rVert \leq \sqrt{\frac{y_i(0)}{\theta_i\lambda_\textrm{max}(R_i)P^2}}\\
&\times \textcolor{black}{\exp\{-\frac{1}{2}\lambda_it\}}, \forall t \in [t_i^h,r]\biggr\}.
\end{align}}

According to (\ref{self3}), we can calculate the value of $\big\lVert \Delta_i(t) \big\rVert$ without using the continuous information. Therefore,  the continuous communication is avoided.  

\textit{Remark 1:}  Since the self-triggered condition (\ref{self4}) is a sufficient condition for the dynamic event-triggered condition (\ref{dynamic ETC}), the number of triggered times by using (\ref{self4})  will be higher than  using (\ref{dynamic ETC}). Our original intention of introducing the self-triggered mechanism is to reduce the consumption of communication resources, which will inevitably increase the number of triggered times.
That is to say, we can  exchange a large amount of communication resources with a small amount of computing resources.

\textit{Remark 2:}  
{\color{black}Compared with the time-triggered methods in \cite{K.G.Vamvoudakis2012}--\cite{H.Zhang2019}, the event-triggered mechanism significantly saves computing resources and communication resources. Note that the event-triggered attitude stabilization problem is studied based on the sliding mode control in \cite{Y.Liu2020}, however, the performance cost has not been consided in the controller design.  } 

\section{Main Results}\label{section4}
Up to now, we have already derived the form of the  optimal controller (\ref{event optimal controller}), which contains the optimal value function $\hat{V}_i^*$.  However, it is very difficult to obtain the analytic solutions to the event-triggered HJB  equations (\ref{event HJB}). 
\textcolor{black}{
In this section, we first introduce an event-triggered Reinforcement Learning (RL) algorithm to obtain the optimal policy. In order to implement the event-triggered RL algorithm online, a  critic neural network is used to approximate the optimal value function $\hat{V}_i^*$. Only measurement data at the
event-triggered instants are needed in the event-triggered
RL algorithm, which obviously  reduces the computation burden.
}

\subsection{Model-Free Event-Triggered RL Algorithm}
\textcolor{black}{
This section presents a model-free event-triggered algorithm based on reinforcement learning, which is used to seek the optimal policy. The RL algorithm involves two parts: policy evaluation and policy improvement. By repeating these two steps at the  event-triggered instants, we know that the optimal policy is obtained when the policy improvement does not change the control policy.
}

\begin{algorithm}\label{event-triggered PI}    
	\caption{Model-Free Event-Triggered RL Algorithm.}   
	Initialize the event-triggered admissible controllers $u_i(0)=0, i=1, ..., N$ and set $h=0$;\\
	\For{\textcolor{black}{the rigid body $i\in \mathcal{V}$} }{   
		\If{\textcolor{black}{the rigid body $i$} receives information $u_j(t)$ transmitted by \textcolor{black}{the rigid body $j$}, where $j \in \mathcal{N}_i$}{
			Set $h'=h'+1$ and $t_j^{h'}=t$;\\
			Update $u_j(t_j^{h'})=u_j(t)$;
		}
		\Else{
			$\hat{u}_j^{h'}=u_j(t_j^{h'})$ remains unchanged;
		}
		Calculate the self-triggered measurment error $\big\lVert \Delta_i (t)\big\rVert$;\\
		\If{the self-triggered condition (\ref{self4}) is violated}{
			Set $t_i^{h+1}=t$;\\
			\textbf{Step 1 (Policy evaluation):}
			\begin{align} \label{algorithm 1-1}
			\nonumber H_i&(e_i,\nabla \hat{V}_i^{h+1},\hat{u}_i^h,\hat{u}_{-i}^{h'})\\
			\nonumber=&e_i^{\top}Q_ie_i+(\hat{u}_i^h)^{\top}R_i\hat{u}_i^h+\big(\nabla \hat{V}_i^{(h+1)}\big)^{\top}\big(X_i\\
			&+l_{ii}Y_i\hat{u}_i^h-\sum_{j\in \mathcal{N}_i}a_{ij}Y_j\hat{u}_j^{h'}\big)=0,
			\end{align}
			where $\nabla\hat{V}_i^{h+1}=\nabla V_i(t_i^{h+1})$, $\hat{u}_i^h=u_i(t_i^h)$ and $\hat{u}_{-i}^{h'}=\hat{u}_j^{h'}= u_{j}(t_j^{h'})$, $j \in \mathcal{N}_i$;\\
			\textbf{Step 2 (Policy Improvement):}
			\begin{align} \label{algorithm 1-2}
			\hat{u}_i^{h+1}=-\frac{1}{2} l_{ii}R_i^{-1} Y_i^{\top} \nabla\hat{V}_i^{h+1},
			\end{align}
			where $\hat{u}_i^{h+1}=u_i(t_i^{h+1})$;\\
			Set $h=h+1$;
		}
		\Else{
			$\hat{u}_i^h=u_i(t_i^h)$ remains unchanged;
		}
	}  
\end{algorithm}

Next, we give a theorem to show  the convergence of the model-free event-triggered RL algorithm.

\textit{Theorem 2:} Supposed that agent $i$ updates the control policy according to Algorithm 1, the value function converges to the optimal value function, i.e., $\lim_{h \to \infty}\hat{V}_i^h=\hat{V}_i^*$ and the control policy converges to the optimal control policy, i.e., $\lim_{h \to \infty} \hat{u}_i^h=\hat{u}_i^*$.

\textit{Proof}: According to Eq. (\ref{algorithm 1-1}), we can obtain that
\begin{align} \label{proof2-1}
\nonumber &\big(\nabla\hat{V}_i^h\big)^{\top} \Big[X_i+l_{ii}Y_i\hat{u}_i^{h-1}-\sum_{j\in \mathcal{N}_i}a_{ij}Y_j\hat{u}_j^{h'-1}\Big]\\
&=-e_i^{\top}Q_ie_i-\big(\hat{u}_i^{h-1}\big)^{\top}R_i\hat{u}_i^{h-1},
\end{align}
and
\begin{align}\label{proof2-2}
\nonumber &\big(\nabla\hat{V}_i^{h+1}\big)^{\top} \Big[X_i+l_{ii}Y_i\hat{u}_i^h-\sum_{j\in \mathcal{N}_i}a_{ij}Y_j\hat{u}_j^{h'}\Big]\\
&=-e_i^{\top}Q_ie_i-\big(\hat{u}_i^h\big)^{\top}R_i\hat{u}_i^{h'}.
\end{align}

By applying the transformation to Eq. (\ref{proof2-1}), the following equation holds:
\begin{align}\label{proof2-3}
\nonumber &\big(\nabla\hat{V}_i^h\big)^{\top} \Big[X_i+l_{ii}Y_i\hat{u}_i^h-\sum_{j\in \mathcal{N}_i}a_{ij}Y_j\hat{u}_j^{h'}\Big]\\
\nonumber&=-e_i^{\top}Q_ie_i-\big(\hat{u}_i^{h-1}\big)^{\top}R_i\hat{u}_i^{h-1}+\big(\nabla\hat{V}_i^h\big)^{\top}l_{ii}Y_i\big(\hat{u}_i^h-\hat{u}_i^{h-1}\big)\\
&\quad-\big(\nabla\hat{V}_i^h\big)^{\top}\sum_{j\in \mathcal{N}_i}a_{ij}Y_j\big(\hat{u}_j^{h'}-\hat{u}_j^{h'-1}\big).
\end{align}

{\color{black}In the model-free event-triggered RL algorithm, we update the control policy of \textcolor{black}{the rigid body $i$} when the self-triggered condition (\ref{self4}) is violated. Under the distributed asynchronous update pattern, the  control policy of \textcolor{black}{the rigid body $j$} remains  invariant, where $j \in \mathcal{N}_i$.} Considering  the trajectory of the consensus error driven by $u_i^h$, i.e., $\dot{e}_i= X_i+l_{ii}Y_i\hat{u}_i^h-\sum_{j\in \mathcal{N}_i}a_{ij}Y_j\hat{u}_j^{h'}$, we have
\begin{align}\label{proof2-4}
\nonumber &\hat{V}_i^{h+1}-\hat{V}_i^h\\
\nonumber&=\int_{t}^\infty\Big[\big(\nabla\hat{V}_i^h\big)^{\top}\dot{e}_i-\big(\nabla\hat{V}_i^{h+1}\big)^{\top}\dot{e}_i\Big]dv\\
\nonumber&=\int_{t}^\infty\Big[\big(\hat{u}_i^h\big)^{\top}R_i\hat{u}_i^h-\big(\hat{u}_i^{h-1}\big)^{\top}R_i\hat{u}_i^{h-1}\\
&\quad+\big(\nabla\hat{V}_i^h\big)^{\top}l_{ii}Y_i\big(\hat{u}_i^h-\hat{u}_i^{h-1}\big)\Big]dv.
\end{align}

According to Eq. (\ref{algorithm 1-2}), it can be easily obtained that
\begin{align}\label{proof2-5}
\big(\nabla\hat{V}_i^h\big)^{\top}l_{ii}Y_i=-2\big(\hat{u}_i^h\big)^{\top}R_i.
\end{align}

Then, Eq. (\ref{proof2-4}) becomes
\begin{align}\label{proof2-6}
\nonumber &\hat{V}_i^{h+1}-\hat{V}_i^h\\
\nonumber&=\int_{t}^\infty\Big[\big(\hat{u}_i^h\big)^{\top}R_i\hat{u}_i^h-\big(\hat{u}_i^{h-1}\big)^{\top}R_i\hat{u}_i^{h-1}\\
\nonumber&\quad-2\big(\hat{u}_i^h\big)^{\top}R_i\big(\hat{u}_i^h-\hat{u}_i^{h-1}\big)\Big]dv\\
\nonumber&=\int_{t}^\infty\Big[-\big(\hat{u}_i^h-\hat{u}_i^{h-1}\big)^{\top}R_i \big(\hat{u}_i^h-\hat{u}_i^{h-1}\big)\Big]dv\\
&\leq 0.
\end{align}

Therefore, $\hat{V}_i^{h+1}\leq \hat{V}_i^{h}$ is always satisfied. According to the Weierstrass theorem \cite{Z.Jiang2013}, the positive definite value function $\hat{V}_i^{h}$ converges to the optimal value function $\hat{V}_i^*$ with $h \to \infty$. Meanwhile, the control policy $\hat{u}_i^h$ converges to the optimal control policy $\hat{u}_i^*$. The complete proof is given.  $\hfill\blacksquare$ 

\subsection{Implementation of  Event-Triggered PI Algorithm}
In this section, we implement Algorithm 1 by using a critic neural network to approximate the optimal value function $\hat{V}_i^*$.

We first define the following neural network of each agent:
\begin{align} \label{critic1}
\hat{V}_i(e_i)=\hat{W}_{c,i}^{\top}\phi_i(e_i), \ \forall t\in[t_i^h,t_i^{h+1}),
\end{align}
where $\hat{W}_{c,i}$ indicates the  critic estimated weight at the event-triggered instant $t_i^h$,   and  $\phi_i(e_i) $ indicates the critic activation function.

According to Eq. (\ref{event HJB}),  the estimated error  of the critic NN can be defined as
\begin{align}\label{critic2}
e_{c,i}=e_iQ_ie_i+\hat{u}_i^{\top}R_i\hat{u}_i+\hat{W}_{c,i}^{\top} \nabla \phi_i \dot{e}_i,
\end{align}
where $\nabla \phi_i=\partial \phi_i(e_i)/\partial e_i^{\top}$.
\begin{figure}[b]
	\centering
	\includegraphics[scale=0.4]{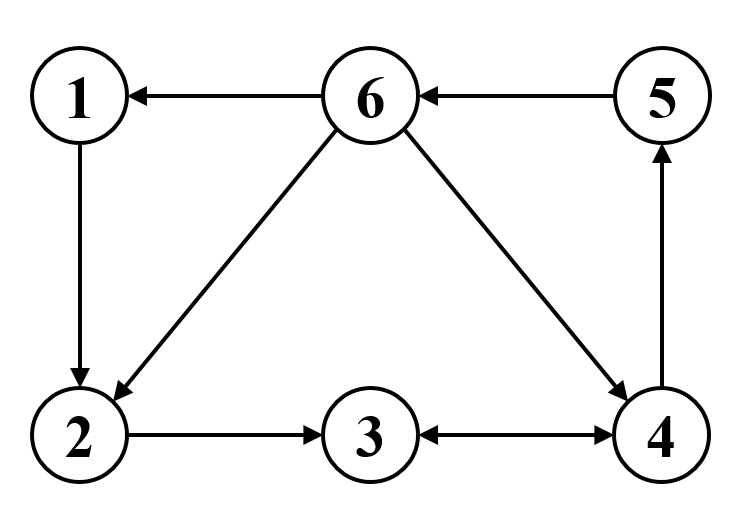}
	\caption{\textcolor{blue}{A strongly connected graph with six nodes. }}
	\label{Communication graph}
\end{figure}
For a given event-triggered admissible controller $\hat{u}_i$, the update rule of $\hat{W}_{c,i}$ is to minimize the following objective function:
\begin{align}\label{critic3}
E_{c,i}=\frac{1}{2}e_{c,i}^{\top}e_{c,i}.
\end{align}

According to the gradient descent method, we can derive the update law of the following form:
\begin{subequations}\label{critic4}
	\begin{align}
	&\dot{\hat{W}}_{c,i}=0,\  t\in (t_i^h,t_i^{h+1}),\\
	&\hat{W}_{c,i}^+=\hat{W}_{c,i}-l_{c,i}k_i(k_{1,i}^{\top}\hat{W}_{c,i}+e_i^{\top}Q_ie_i+\hat{u}_i^{\top}R_i\hat{u}_i), \ t=t_i^h,
	\end{align}	
\end{subequations}
where $l_{c,i}>0$ indicates the learning rate of the critic NN, $k_{1,i}=\nabla
\phi_i\dot{e}_i$ and $k_i=k_{1,i}/(k_{1,i}^{\top} k_{1,i}+1)^2$.

Letting $\tilde{W}_{c,i}=\hat{W}_{c,i}-W_{c,i}$,  we can deduce
\begin{subequations}\label{critic5}
	\begin{align}
	&\dot{\tilde{W}}_{c,i}=0,\  t\in (t_i^h,t_i^{h+1}),\\
	&\tilde{W}_{c,i}^+=\tilde{W}_{c,i}-l_{c,i}k_i(k_{1,i}^{\top}\tilde{W}_{c,i}+\epsilon_{c,i}), \  t=t_i^h,
	\end{align}	
\end{subequations}
where $W_{c,i}$ denotes the  critic target weight, $\tilde{W}_{c,i}$ is the critic weight error and $\epsilon_{c,i}=e_iQ_ie_i+\hat{u}_i^{\top}R_i\hat{u}_i+W_{c,i}^{\top} \nabla \phi_i \dot{e}_i$ is the critic residual error.

Therefore, the optimal controller can be obtained by (\ref{algorithm 1-2}) and (\ref{critic1}), which is expressed in the following form:
\begin{align}
\hat{u}_i=-\frac{1}{2} l_{ii}R_i^{-1} Y_i^{\top} (\nabla \phi_i)^{\top} \hat{W}_{c,i}.
\end{align}

Through the above critic NN framework, we can obtain the optimal controller with only measurement data at the event-triggered instants. Therefore, the need for system dynamics is obviously avoided.  In addition, the neural network is only updated at the event-triggered instants $t_i^h$,  which are determined by  the self-triggered condition (\ref{self4}).

\textit{Assumption 3:} In the critic NN framework, the target weight matrix,  the activation function, the critic residual error are bounded with positive constants $W_{cM}$, $\phi_M$, and $\epsilon_{cM}$, i.e., 
$\lVert W_{c,i} \rVert \leq W_{cM}$,
$\lVert \phi_i \rVert \leq \phi_M$, and
$\lVert \epsilon_{c,i} \rVert \leq \epsilon_{cM}$.

\textit{Theorem 3:} Consider the consensus error dynamics (\ref{event dot consensus error}), the critic neural network is given as (\ref{critic1}). If the estimated weight matrix $\hat{W}_{c,i}$ is updated with (\ref{critic4}), the consensus error $e_i$  and the critic estimation error $\tilde{W}_{c,i}$ are UUB.

\begin{figure*}[htb]
	\centering
	\subfigure[]{
		\includegraphics[width=3.2in]{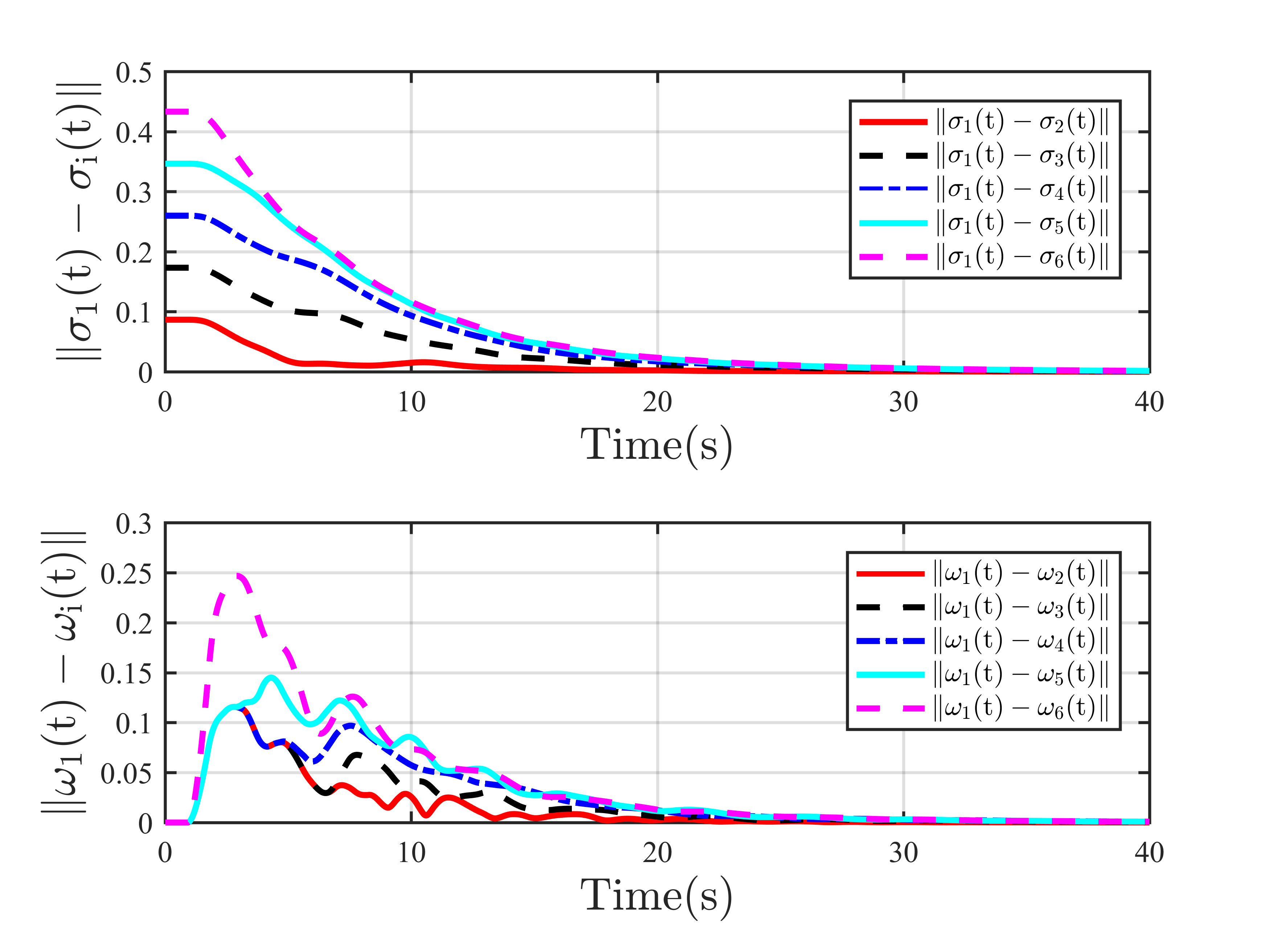}
		\label{attitude}
	}
	\hfill
	\subfigure[]{
		\includegraphics[width=3.2in]{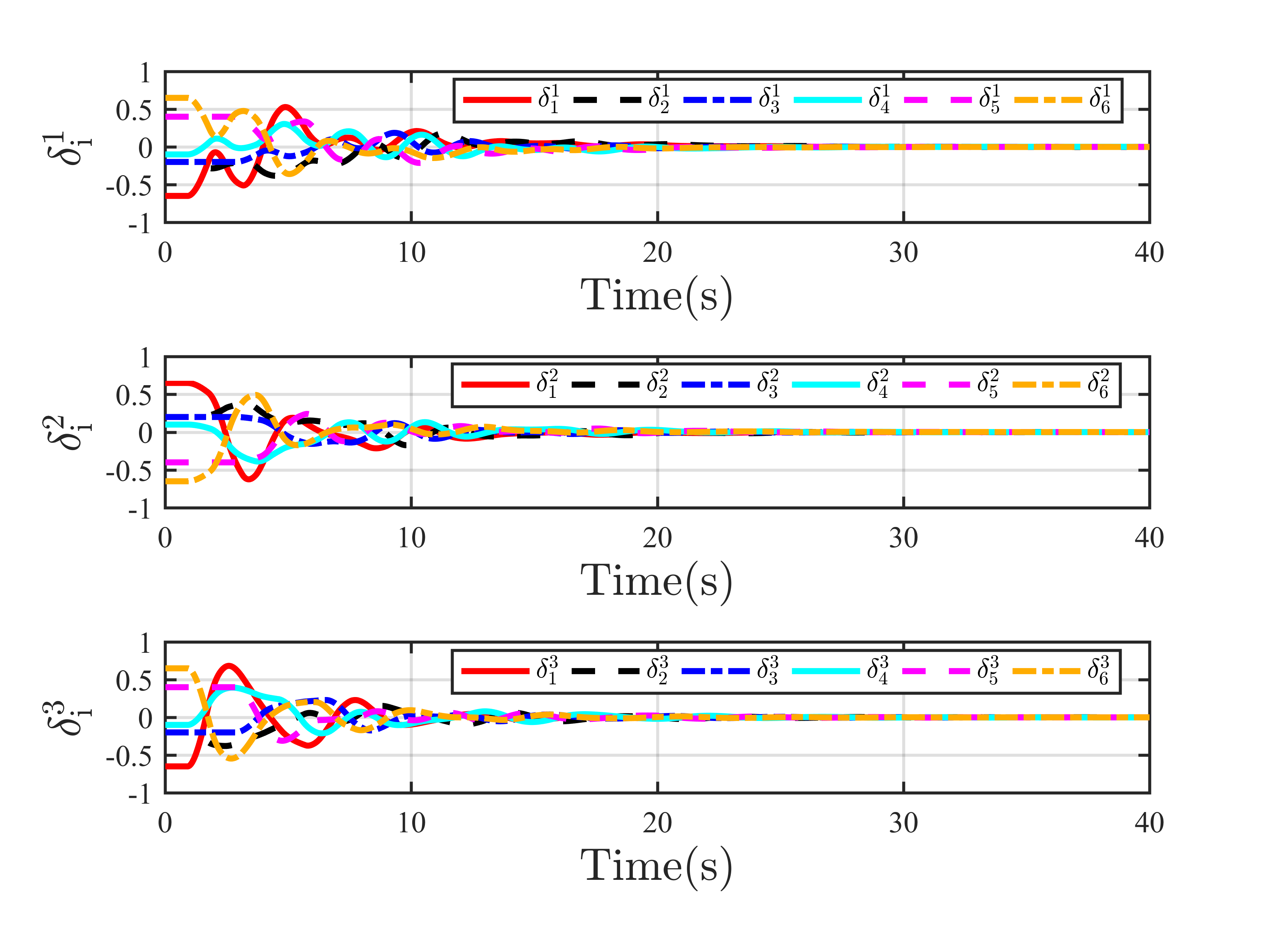}
		\label{delta}
	}
	\caption{\textcolor{blue}{(a) The norms of attitude errors and angular velocity errors. The trajectories indicate 
		the norms of attitude errors between the rigid body $1$ and the rigid body $i \in \{2,3,4,5,6\}$.
		(b) The consensus errors of each rigid body $\delta_i$. Three subfigures show three components of the consensus error vector $\delta_i$ of each agent, respectively. }} 
	\end{figure*}
\begin{figure*}[htb]
	\subfigure[]{
		\includegraphics[width=3.2in]{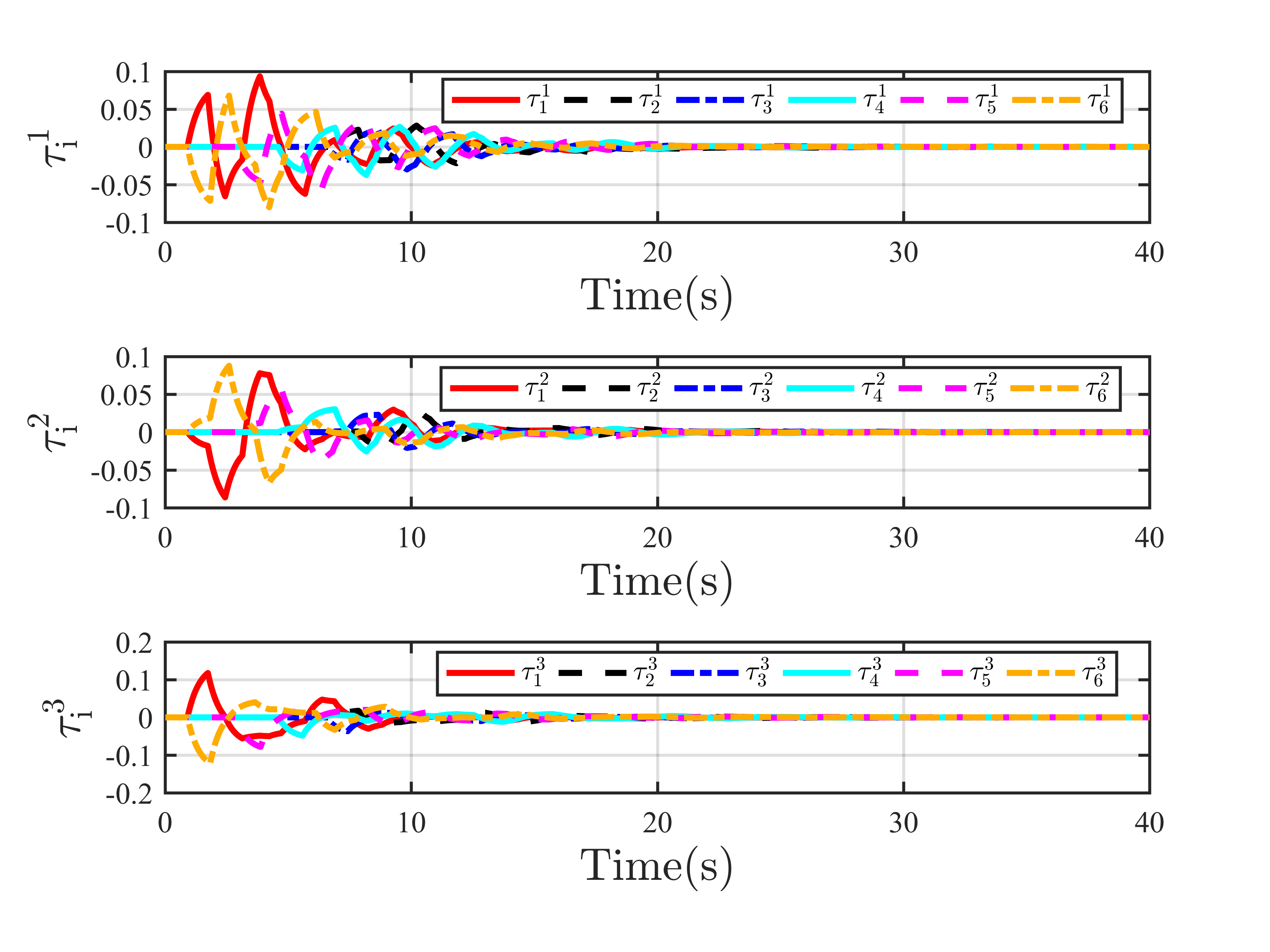}
		\label{tau}
	}
	\hfill
	\subfigure[]{
		\includegraphics[width=3.2in]{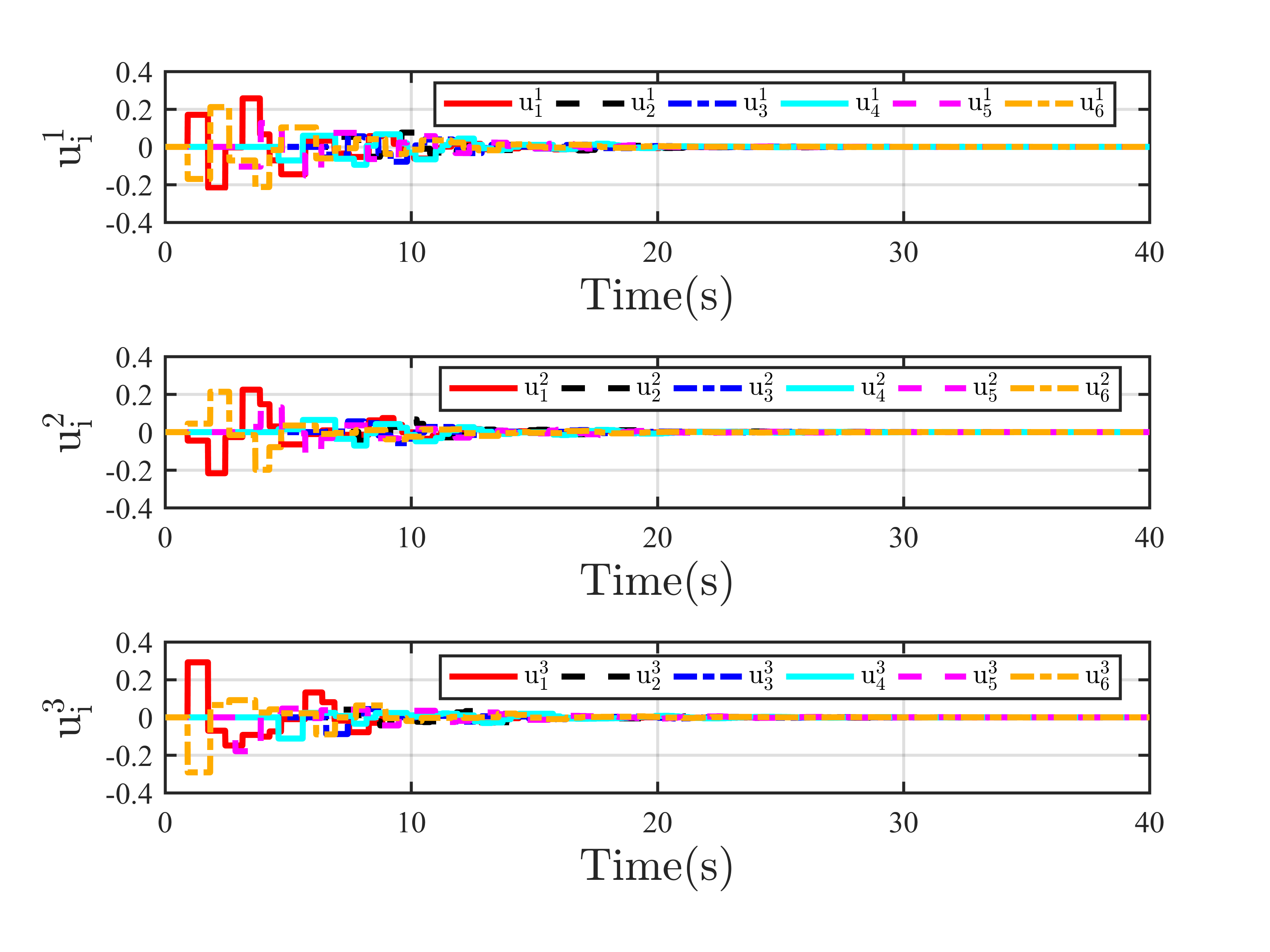}
		\label{u}
	}
	\caption{ \textcolor{blue}{(a) The original control inputs $\tau_i$ of each rigid body.  Three subfigures show three components of the control input $\tau_i$ of each agent, respectively. 
		(b) The control inputs $\text{u}_i $ of the augmented systems. Three subfigures show three components of the control input $\text{u}_i $ of augmented systems of each agent, respectively. } }
\end{figure*}
\begin{figure}[b]
	\centering
	{
		\includegraphics[scale=0.5]{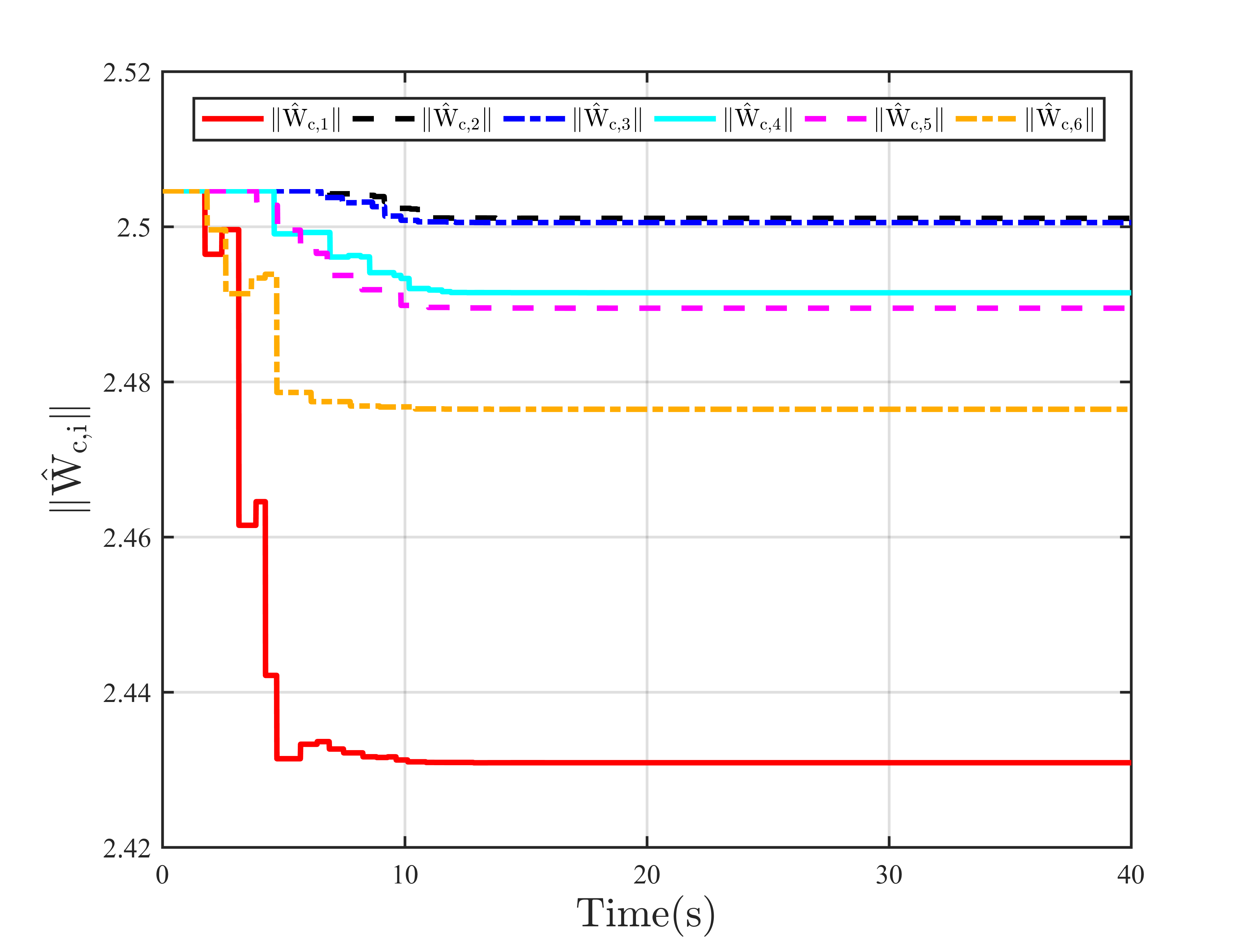}
		\caption{\textcolor{blue}{The critic estimated weight matrices. The trajectories show the norm of weight matrices of each agent.}}
		\label{W}
	}
\end{figure}

\textit{Proof}: Two different situations are considered, including  during the event-triggered intervals and at the event-triggered instants. 

\emph{Situation 1:} During the event-triggered intervals, i.e., $t \in (t_i^h, t_i^{h+1})$.

Consider the Lyapunov function of the following form:
\begin{align}\label{proof3-1}
L_i=L_{i,1}+L_{i,2},
\end{align}
where $L_{i,1}=e_i^{\top}e_i+V_i(e_i)$, $L_{i,2}=\frac{\text{tr}(\tilde{W}_{c,i}^{\top}\tilde{W}_{c,i})}{l_{c,i}}$.

According to (\ref{critic5}), we can obtain
\begin{align}\label{proof3-2}
\dot{L}_{i,2}=\frac{2\text{tr}(\tilde{W}_{c,i}^{\top}\dot{\tilde{W}}_{c,i})}{l_{c,i}}=0.
\end{align}

Therefore, the first-order derivative of $L_i$ can be expressed as follows:
\begin{align}\label{proof3-3}
\nonumber\dot{L}_i&=\dot{L}_{i,1}=2e_i^{\top}\dot{e}_i+\dot{V}_i(e_i)\\
\nonumber&=2e_i^{\top}\big(X_i(e_i)+l_{ii}Y_i\hat{u}_i-\sum_{j\in \mathcal{N}_i}a_{ij}Y_j\hat{u}_j\big)\\
\nonumber&\quad-e_i^{\top}Q_ie_i-\hat{u}_i^{\top}R_i\hat{u}_i\\
\nonumber&\leq\lVert e_i \rVert^2+\big\lVert X_i(e_i)+l_{ii}Y_i\hat{u}_i-\sum_{j\in \mathcal{N}_i}a_{ij}Y_j\hat{u}_j \big\rVert^2\\
\nonumber&\quad-\lambda_\textrm{min}(Q_i)\lVert e_i \rVert^2-\lambda_\textrm{min}(R_i)\lVert \hat{u}_i \rVert^2\\
\nonumber&\leq\big(1+3X_M^2-\lambda_\textrm{min}(Q_i)\big)\lVert e_i \rVert^2+3l_{ii}^2Y_M^2 \lVert \hat{u}_i \rVert^2\\
&\quad+3\sum_{j\in \mathcal{N}_i}a_{ij}^2Y_M^2\lVert \hat{u}_j \rVert^2-\lambda_\textrm{min}(R_i)\lVert \hat{u}_i \rVert^2.
\end{align}

In order to ensure $\dot{L}_i < 0$, the following inequality should be satisfied:
\begin{align}\label{proof3-4}
\lVert e_i \rVert>\sqrt \frac{\Phi_i}{\lambda_\textrm{min}(Q_i)-1-3X_M^2},
\end{align}
where $\Phi_i=3l_{ii}^2Y_M^2 \lVert \hat{u}_i \rVert^2+3\sum_{j\in \mathcal{N}_i}a_{ij}^2Y_M^2\lVert \hat{u}_j \rVert^2-\lambda_\textrm{min}(R_i)\lVert \hat{u}_i \rVert^2$.

Hence, the consensus error $e_i$ is UUB. During the event-triggered intervals, the critic estimation error $\tilde{W}_{c,i}$ remains unchanged, which means $\tilde{W}_{c,i}$ is also UUB.

\emph{Situation 2:} At the event-triggered instants, i.e., $t=t_i^h$.

Choosing the same Lyapunov function as (\ref{proof3-1}), we can obtain:
\begin{align}\label{proof4-1}
\Delta L_i=\Delta L_{1,i}+\Delta L_{2,i}.
\end{align}

Since the trajectory of $e_i$ is continuous, i.e., $e_i^+=e_i$. Therefore, one has
\begin{align}\label{proof4-2}
\Delta L_{1,i}=(e_i^+)^{\top}e_i^++V_i(e_i^+)-e_i^{\top}e_i-V_i(e_i)=0.
\end{align}

Next, according to (\ref{critic5}), we have
\begin{align}\label{proof4-3}
\nonumber\Delta L_{2,i}&=\frac{\text{tr}\big[(\tilde{W}_{c,i}^+)^{\top}\tilde{W}_{c,i}^+\big]}{l_{c,i}}-\frac{\text{tr}(\tilde{W}_{c,i}^{\top}\tilde{W}_{c,i})}{l_{c,i}}\\
\nonumber&=\frac{1}{l_{c,i}}\text{tr}\Big[\big(\tilde{W}_{c,i}-l_{c,i}k_i(k_{1,i}^{\top}\tilde{W}_{c,i}+\epsilon_{c,i})\big)^{\top}\\
\nonumber&\quad\times\big(\tilde{W}_{c,i}-l_{c,i}k_i(k_{1,i}^{\top}\tilde{W}_{c,i}+\epsilon_{c,i})\big)- \tilde{W}_{c,i}^{\top}\tilde{W}_{c,i}\Big]\\
\nonumber&=l_{c,i}\text{tr}\Big[\big(k_i(k_{1,i}^{\top}\tilde{W}_{c,i}+\epsilon_{c,i})\big)^{\top}\big(k_i(k_{1,i}^{\top}\tilde{W}_{c,i}+\epsilon_{c,i})\big)\Big]\\
\nonumber&\quad-2\text{tr}(\tilde{W}_{c,i}^{\top}k_ik_{1,i}^{\top}\tilde{W}_{c,i})-2\text{tr}(\tilde{W}_{c,i}^{\top}k_i\epsilon_{c,i})\\
\nonumber&=l_{c,i}\lVert k_ik_{1,i}^{\top}\tilde{W}_{c,i}+k_i\epsilon_{c,i}\rVert^2-2k_i^{\top}k_{1,i}\lVert \tilde{W}_{c,i}\rVert^2\\
&\quad-2\lVert\tilde{W}_{c,i}^{\top}k_i\epsilon_{c,i}\rVert.
\end{align}

From the definition of $k_{1,i}$ and $k_i$,  we can obtain the following  inequalities:
\begin{align}\label{proof4-4}
\alpha_k \leq k_{1,i}^{\top}&k_i \leq \beta_k,\\
\lVert k_i\rVert \leq & K_M,
\end{align}
where $\beta_k>\alpha_k>0$ and $K_M>0$.

Substituting (\ref{proof4-4}) and  (61) into (\ref{proof4-3}), $\Delta L_{2,i}$ becomes
\begin{align}\label{proof4-5}
\nonumber\Delta L_{2,i} &\leq 2 l_{c,i}\beta_k^2\lVert\tilde{W}_{c,i}\rVert^2+2l_{c,i}\epsilon_{cM}^2K_M^2-2\alpha_k\lVert\tilde{W}_{c,i}\rVert^2\\
\nonumber&\quad+\epsilon_{cM}(\lVert\tilde{W}_{c,i}\rVert^2+K_M^2)\\
\nonumber&\leq-(2\alpha_k-2 l_{c,i}\beta_k^2-\epsilon_{cM})\lVert\tilde{W}_{c,i}\rVert^2\\
&\quad+(2l_{c,i}\epsilon_{cM}^2+\epsilon_{cM})K_M^2.
\end{align}

Combining (\ref{proof4-2}) and (\ref{proof4-5}), $\Delta L_i$ can be transformed into the following form:
\begin{align}\label{proof4-6}
\nonumber\Delta L_i&=\Delta L_{1,i}+\Delta L_{2,i}\\
\nonumber&\leq-(2\alpha_k-2 l_{c,i}\beta_k^2-\epsilon_{cM})\lVert\tilde{W}_{c,i}\rVert^2\\
&\quad+(2l_{c,i}\epsilon_{cM}^2+\epsilon_{cM})K_M^2.
\end{align}

In order to simplify  the expression, some auxiliary variables are defined as follows:
\begin{align}
\nonumber A_i&=2\alpha_k-2 l_{c,i}\beta_k^2-\epsilon_{cM},\\
\nonumber \Gamma_i&=(2l_{c,i}\epsilon_{cM}^2+\epsilon_{cM})K_M^2.
\end{align}

Therefore, it can be deduced $\Delta L_i < 0$ when $\lVert \tilde{W}_{c,i} \rVert > \sqrt{\Gamma_i/A_i}$, which signifies that $e_i$ and $\tilde{W}_{c,i}$ are UUB at the event-triggered instants.

Combing situation 1 with situation 2,  it can be proved that the consensus error $e_i$ and the critic estimation error $\tilde{W}_{c,i}$ are UUB. The complete proof is given. $\hfill\blacksquare$

{\color{black}\textit{Remark 3:} The attitude consensus problem of multiple rigid body networks has been widely studied in the literature \cite{A.Abdessameud2009,H.Cai2016,H.Gui2018,M.Lu2020}. However, the attitude consensus protocol proposed in \cite{A.Abdessameud2009,H.Cai2016,H.Gui2018,M.Lu2020} are all based on the known rigid body dynamics, which is a major limitation in practical applications. 
In this work, a model-free RL algorithm is proposed to solve the HJB equation of the optimal attitude consensus of multiple rigid body networks.
Moreover, compared with the existing results on model-free consensus problem of multi-agent networks \cite{J.Li2017, J.Qin2019, H.Zhang2019}, an event-triggered RL algorithm is proposed, which is further extended to the self-triggered RL algorithm. 
Based on Algorithm 1, we know that the control update action and the information interaction among agents are only executed on the triggering instants. Hence, the computation and communication resource can be obviously reduced compared with the continuous-time approaches \cite{J.Li2017, J.Qin2019, H.Zhang2019}. }


\section{Simulation}\label{section5}
This section presents a numerical simulation to verify the effectiveness of the proposed event-triggered reinforcement learning method.  We consider a multiple rigid body network with six nodes under the strongly connected communication topology. The communication relationship between any two nodes can be seen in Fig. \ref{Communication graph}. The Laplacian matrix is selected as follows:
\begin{align}
\nonumber \mathcal{L}=\begin{bmatrix}4&0&0&0&0&-4\\-4&8&0&0&0&-4\\0&-4&8&-4&0&0\\0&0&-4 &8&0&-4\\0&0&0&-4&4&0\\0&0&0&0&-4&4\end{bmatrix}.
\end{align}
\begin{figure*}
	\centering
	\subfigure[]{
		\includegraphics[width=3.2in]{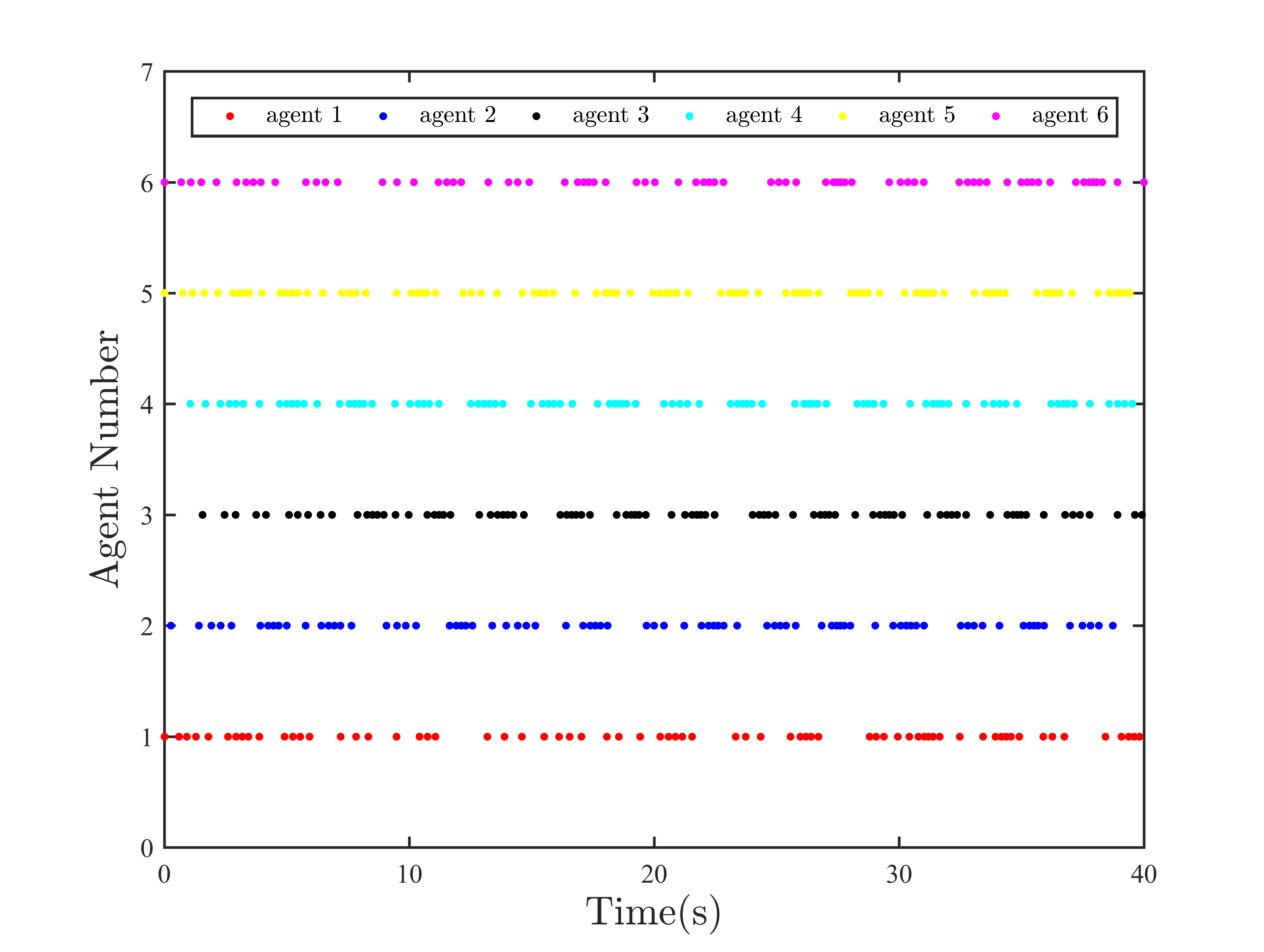}
		\label{Triggering_instants1}
	}
	\hfill
	\subfigure[]{
		\includegraphics[width=3.2in]{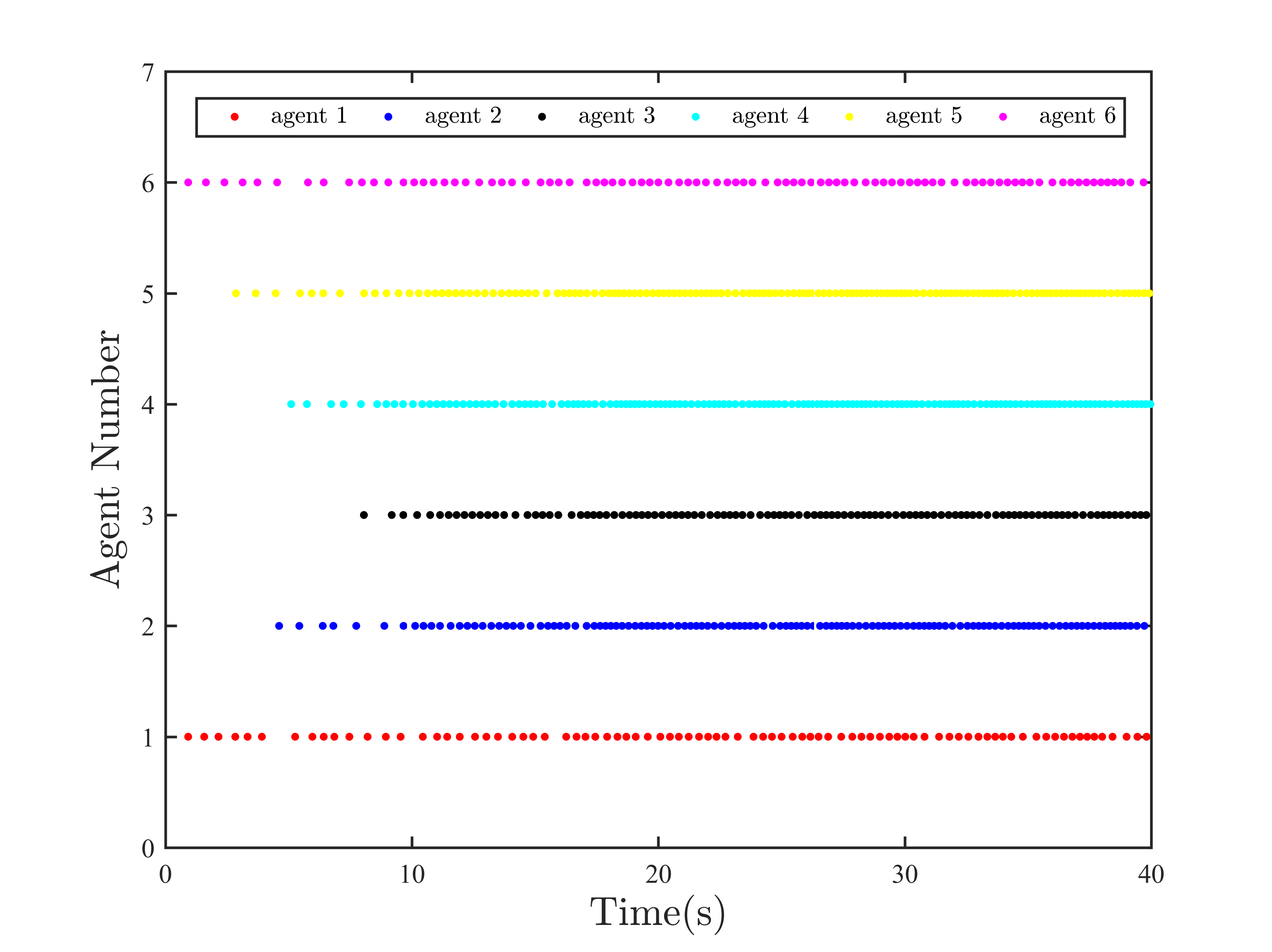}
		\label{Triggering_instants2}
	}
	\caption{ (a) Triggering instants of dynamic event-triggered control of each agent. (b) Triggering instants of self-triggered control of each agent.}
\end{figure*}
\begin{figure*}[htb]
	\centering
	\subfigure[]{
		\includegraphics[width=3.3in]{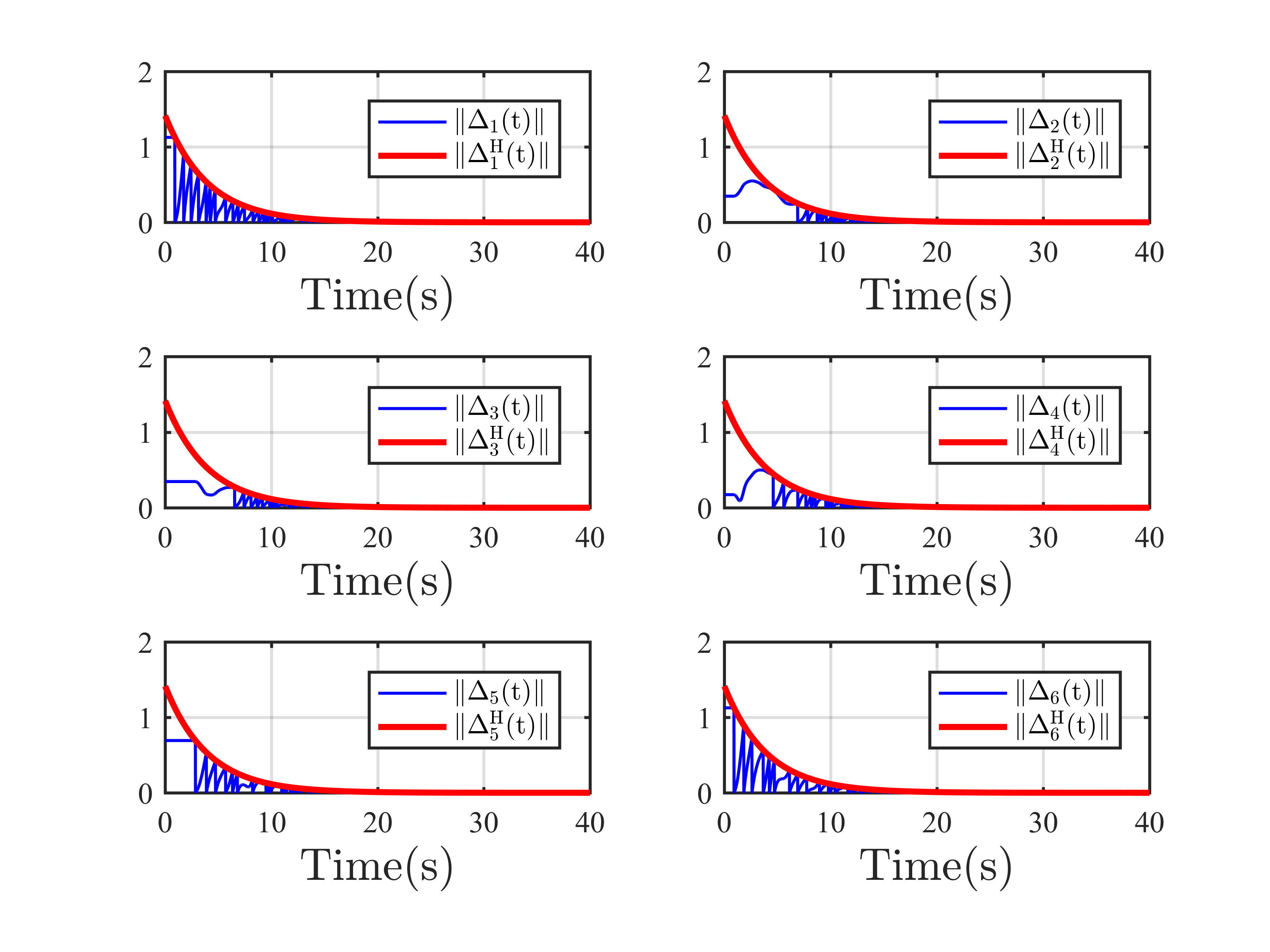}
		\label{Delta}
	}
	\hfill
	\subfigure[]{
		\includegraphics[width=3.2in]{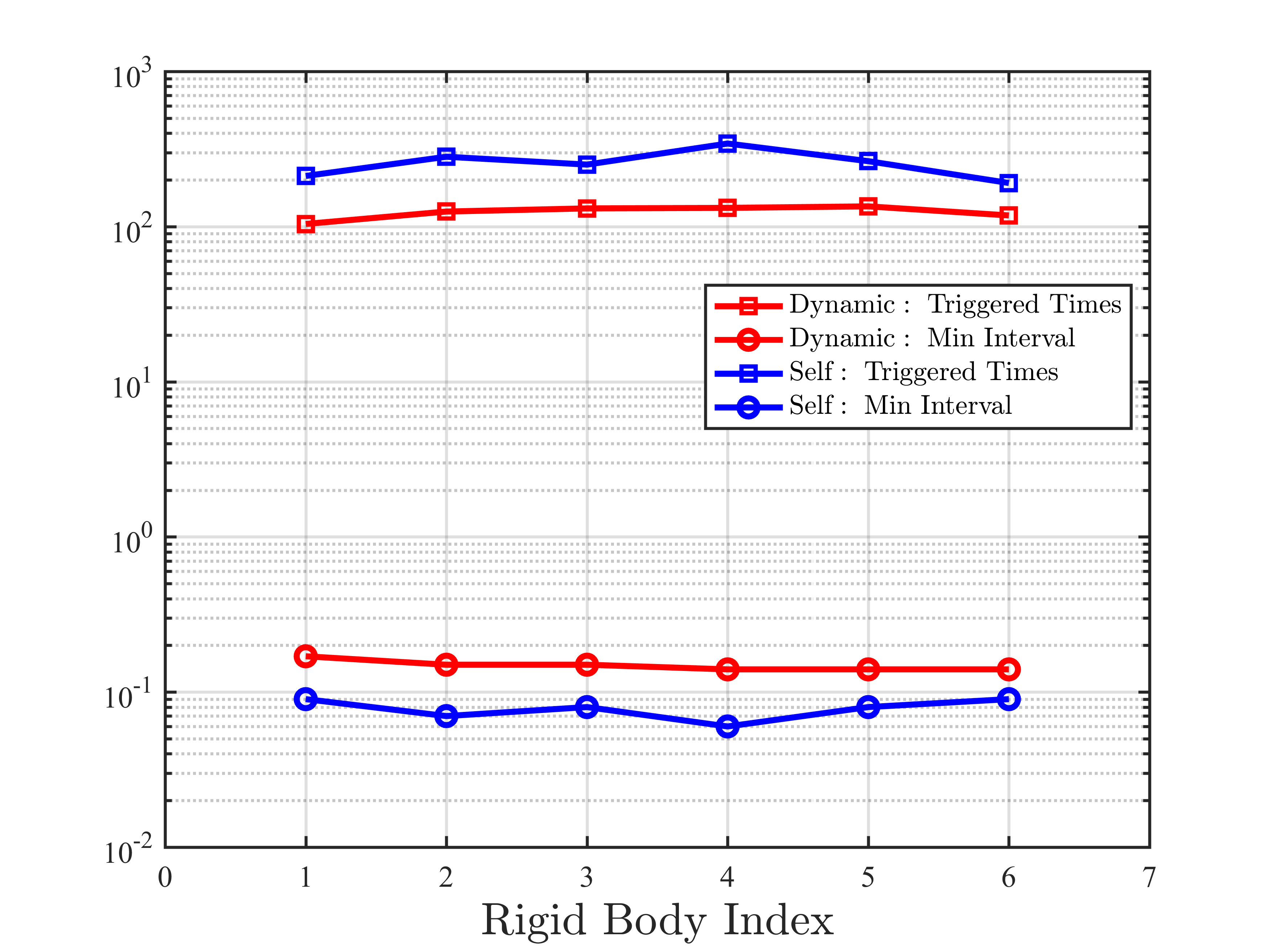}
		\label{compare}
	}
	
	\caption{ (a) The self-triggered measurement errors and their  upper bound of each agent. (b) Comparisons of the total number of triggering times and minimum triggering intervals between the dynamic event-triggered control and the self-triggered control of each agent. }
\end{figure*}
{\color{black}The  rigid body $i \in \{1,2,3,4,5,6\}$ can be modeled using the following equations:
\begin{align}
	\nonumber \dot{\sigma_i}&=G(\sigma_i)\omega_i,\\
	\nonumber  J_i\dot{\omega}_i&=-\omega_i \times (J_i\omega_i)+\tau_i,
\end{align}
in which $\sigma_i=[\sigma_i^{(1)}, \sigma_i^{(2)}, \sigma_i^{(3)}]^\top$ indicates the attitude vector, $\omega_i=[\omega_i^{(1)}, \omega_i^{(2)}, \omega_i^{(3)}]^\top$ indicates the angular velocity vector, $J_i$ indicates the inertial matrix. Here, the inertial matrix of each rigid body is selected as follows:
\begin{align}
	\nonumber &J_1=J_3=[1.0\ 0.1\ 0.1;\ 0.1\ 1.0\ 0.1;\ 0.1\ 0.1\ 1.0],\\
	\nonumber &J_2=J_4=[1.2\ 0.1\ 0.1;\ 0.1\ 0.9\ 0.1;\ 0.1\ 0.1\ 1.1],\\
	\nonumber &J_3=J_6=[1.1\ 0.2\ 0.1;\ 0.2\ 1.0\ 0.3;\ 0.1\ 0.3\ 1.3].
\end{align}

In this simulation, the total duration is set to 40 seconds and the sampling period is 0.01 seconds. We selected the parameters as follows: $\alpha_i = 0.5$, $P = 1$, the weight matrices $Q_i=4I_6$, $R_i=I_3$, and the learning rate $l_{c,i}=0.6$. In the dynamic event-triggered condition (\ref{dynamic ETC}), $y_i(0)=4$, $\gamma_i=0.5$, $\kappa_i=0.5$, $\varpi_i=0.6$, $\theta_i=2$. In the self-triggered condition (\ref{self4}), the parameters remain the same except $\kappa_i=0$ and $\varpi_i=0$. The initial states of each rigid body are given by:
\begin{equation}
	\begin{split}
		\nonumber
		&\sigma_i=\begin{bmatrix}0.05i\\-0.05i\\0.05i\end{bmatrix},\;\omega_i=\dot{\omega}_i=\begin{bmatrix}0\\0\\0\end{bmatrix},\; i=1,2,3,4,5,6.
	\end{split}
\end{equation}}
The critic activation function is designed as:
\begin{align}
\nonumber\phi_i(e_i)=[&(e_i^1)^2 \quad e_i^1e_i^2 \quad e_i^1e_i^3 \quad e_i^1e_i^4 \quad e_i^1e_i^5 \quad e_i^1e_i^6\\
\nonumber &(e_i^2)^2 \quad e_i^2e_i^3 \quad e_i^2e_i^4 \quad e_i^2e_i^5 \quad e_i^2e_i^6 \quad (e_i^3)^2\\
\nonumber &e_i^3e_i^4 \quad e_i^3e_i^5 \quad e_i^3e_i^6 \quad (e_i^4)^2 \quad e_i^4e_i^5 \quad e_i^4e_i^6\\
\nonumber &(e_i^5)^2 \quad e_i^5e_i^6 \quad (e_i^6)^2]^{\top} \in \mathbb{R}^{21}.
\end{align}

By using the model-free event-triggered RL method proposed above, the optimal attitude consensus problem for multiple rigid body networks is solved. Fig. \ref{attitude} shows  the norms of attitude errors and angular velocity errors between \textcolor{black}{the rigid body $1$ and the rigid body $i \in \{2,3,4,5,6\}$}.
From Fig. \ref{attitude}, we can conclude  that the optimal attitude consensus is achieved.
The same conclusion can be obtained from Fig. \ref{delta}, which shows the trajectories of the consensus errors.
Fig. \ref{tau} and Fig. \ref{u} show the original control inputs of each rigid body and the control inputs of the augmented systems, respectively. It is worth noting that the control inputs of the augmented systems are only updated at the event-triggered instants, and remain unchanged during the event-triggered intervals.

Fig. \ref{W} demonstrates the critic estimated weight matrix  of each rigid body. It can be clearly seen that the neural networks are only updated at the event-triggered instants, which obviously reduces the consumption of computing resources.
{\color{black}The triggering instants of dynamic event-triggered control and self-triggered control are illustrated in Figs. \ref{Triggering_instants1} and \ref{Triggering_instants2}, respectively. 
It is clearly shown that control update actions and communication frequencey are both significantly reduced compared with the continous-time control approaches \cite{J.Li2017,J.Qin2019,H.Zhang2019}. 
Fig. \ref{Delta} shows the self-triggered measurement errors and their  upper bound, which determines the event-triggered instants.
Fig. \ref{compare} represents the triggered times and the minimum triggered interval under the dynamic event-triggered condition and the self-triggered condition, respectively. 
Since the self-triggered measurement error $\Delta_i(t)$ is the upper bound of the measurement error $E_i(t)$, the triggered times under the self-triggered mechanism are more than uner the dynamic event-triggered mechanism.  Therefore, we can conclude that the self-triggered mechanism leads to an inevitable increase in the number of triggered times without continuing to communicate with neighbors.}

\section{Conclusion}\label{section6}
In this paper, a model-free event-triggered RL method is proposed to deal with the optimal attitude consensus  for multiple rigid body networks, which only requires the measurement data at the event-triggered instants. In order to solve the HJB equations, an event-triggered PI algorithm is proposed to obtain the optimal policy. Meanwhile, the critic NN framework is used to approximate the optimal value function online. The critic neural network is updated only when the event-triggered condition is violated, which greatly reduces the consumption of computing and communication resources.
The UUB of the consensus error and the weight estimation error is proved and the Zeno behavior is excluded. A numerical simulation for a multiple rigid body network with six nodes  shows the feasibility of the proposed method.

\textcolor{black}{
In the future, we will further improve this work from the following perspectives. 
One consideration is to relax the condition of communication topologies, such as from strongly connected graphs to directed spanning trees or even switching topologies \cite{L.Lei2020}. Due to the actuator failure could happen in real applications of rigid body such as  intelligent cars and
quadrotor aircrafts \cite{Y.Liu2020_TIE}, it is well motivated to consider the optimal
cooperative control of rigid body systems with the actuator failure.
}



\ifCLASSOPTIONcaptionsoff
\newpage
\fi



%

\begin{IEEEbiography}[{\includegraphics[width=1in,height=1.25in,clip,keepaspectratio]{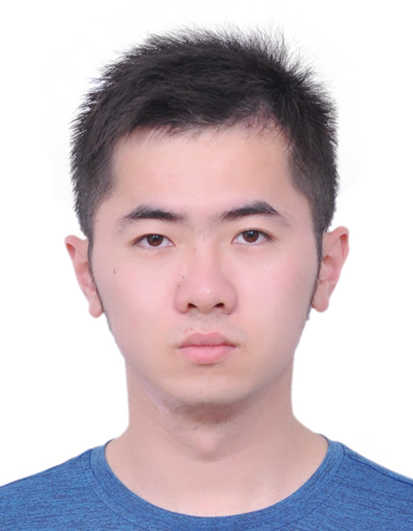}}]{Xin Jin} received the B.S. degree in school of automation from the Guangdong University of Technology, Guangzhou, China, in 2016. He was an exchange Ph.D. student at University of Victoria, Victoria, Canada from Sept. 2019 to Sept. 2020. He is currently working toward the Ph.D. degree from the East China University of Science and Technology. His research interests include rigid body systems, multi-agent systems, event-triggered control and their applications.
\end{IEEEbiography}

\begin{IEEEbiography}[{\includegraphics[width=1in,height=1.25in,clip,keepaspectratio]{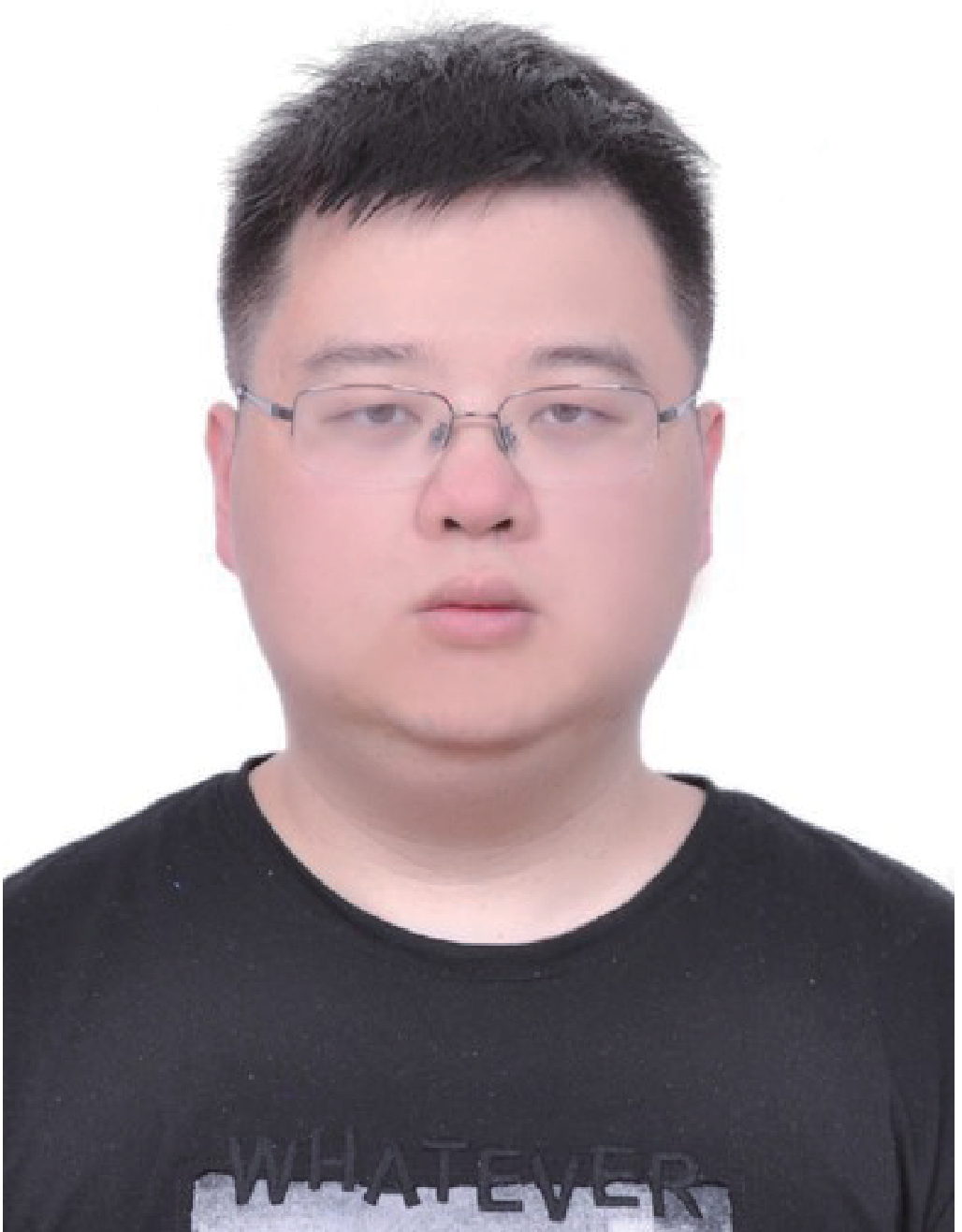}}]{Shuai Mao}
	received the B.S. degree in school of control science and engineering from East China University of Science and Technology, in 2017. He is currently pursuing the Ph.D. degree from East China University of Science and Technology. His research interests include multi-agent systems, distributed optimization and their applications in practical engineering.
\end{IEEEbiography}

\begin{IEEEbiography}[{\includegraphics[width=1in,height=1.25in,clip,keepaspectratio]{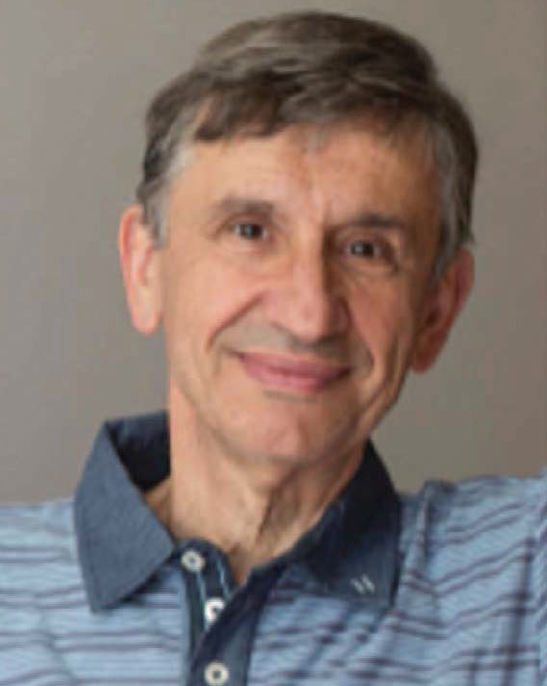}}]{Ljupco Kocarev}
(Fellow, IEEE) is currently a member
of the Macedonian Academy of Sciences and
Arts, a Full Professor with the Faculty of Computer
Science and Engineering, Ss. Cyril and Methodius
University, Skopje, Macedonia, the Director of the
Research Center for Computer Science and Information
Technologies, Macedonian Academy, and a
Research Professor with the University of California
at San Diego. His work has been supported by
the Macedonian Ministry of Education and Science,
the Macedonian Academy of Sciences and Arts,
NSF, AFOSR, DoE, ONR, ONR Global, NIH, STMicroelectronics, NATO,
TEMPUS, FP6, FP7, Horizon 2020, and agencies from Spain, Italy, Germany
(DAAD and DFG), Hong Kong, and Hungary. His scientific interests include
networks, nonlinear systems and circuits, dynamical systems and mathematical
modeling, machine learning, and computational biology.
\end{IEEEbiography}
\begin{IEEEbiography}[{\includegraphics[width=1in,height=1.25in,clip,keepaspectratio]{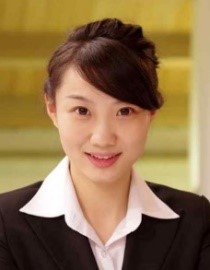}}]{Chen Liang}
	is currently working as a Research Assistant at Key Laboratory of Smart Manufacturing in Energy Chemical Process, Ministry of Education and a Faculty Member of School of Information at East China University of Science and Technology. She got her Master Degree in Computer Applied Technology at Shanghai Normal University in 2013. Her research interests include multi-agent systems, reinforcement learning and network.
	\end{IEEEbiography}
\begin{IEEEbiography}[{\includegraphics[width=1in,height=1.25in,clip,keepaspectratio]{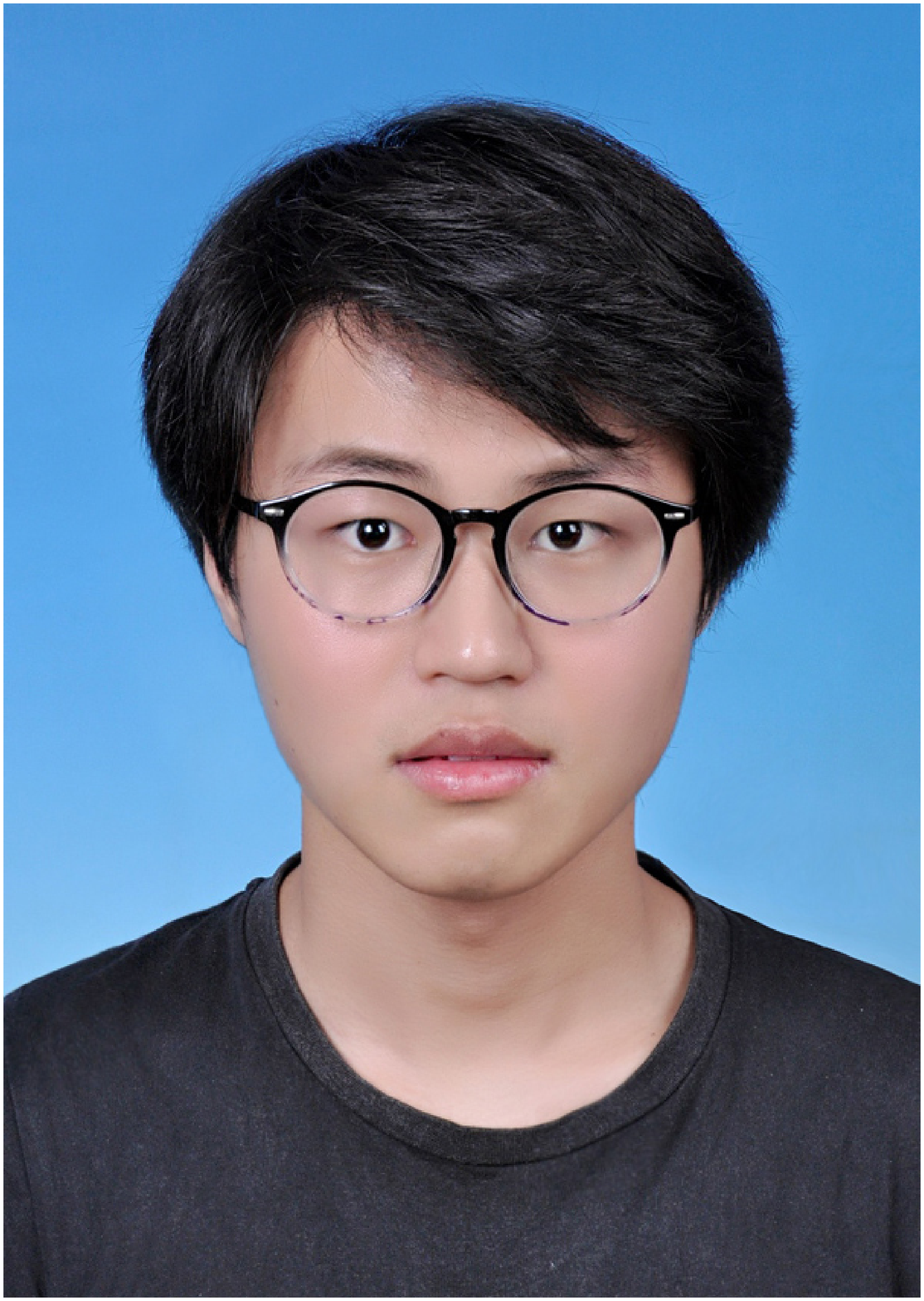}}]{Saiwei Wang}
	received the B.S. degree in school of control science and engineering from East China University of Science and Technology, Shanghai, China, in 2018. He is currently pursuing the M.S. degree from East China University of Science and Technology. His research interests include multi-agent systems, reinforcement learning and their applications.
\end{IEEEbiography}

\begin{IEEEbiography}[{\includegraphics[width=1in,height=1.25in,clip,keepaspectratio]{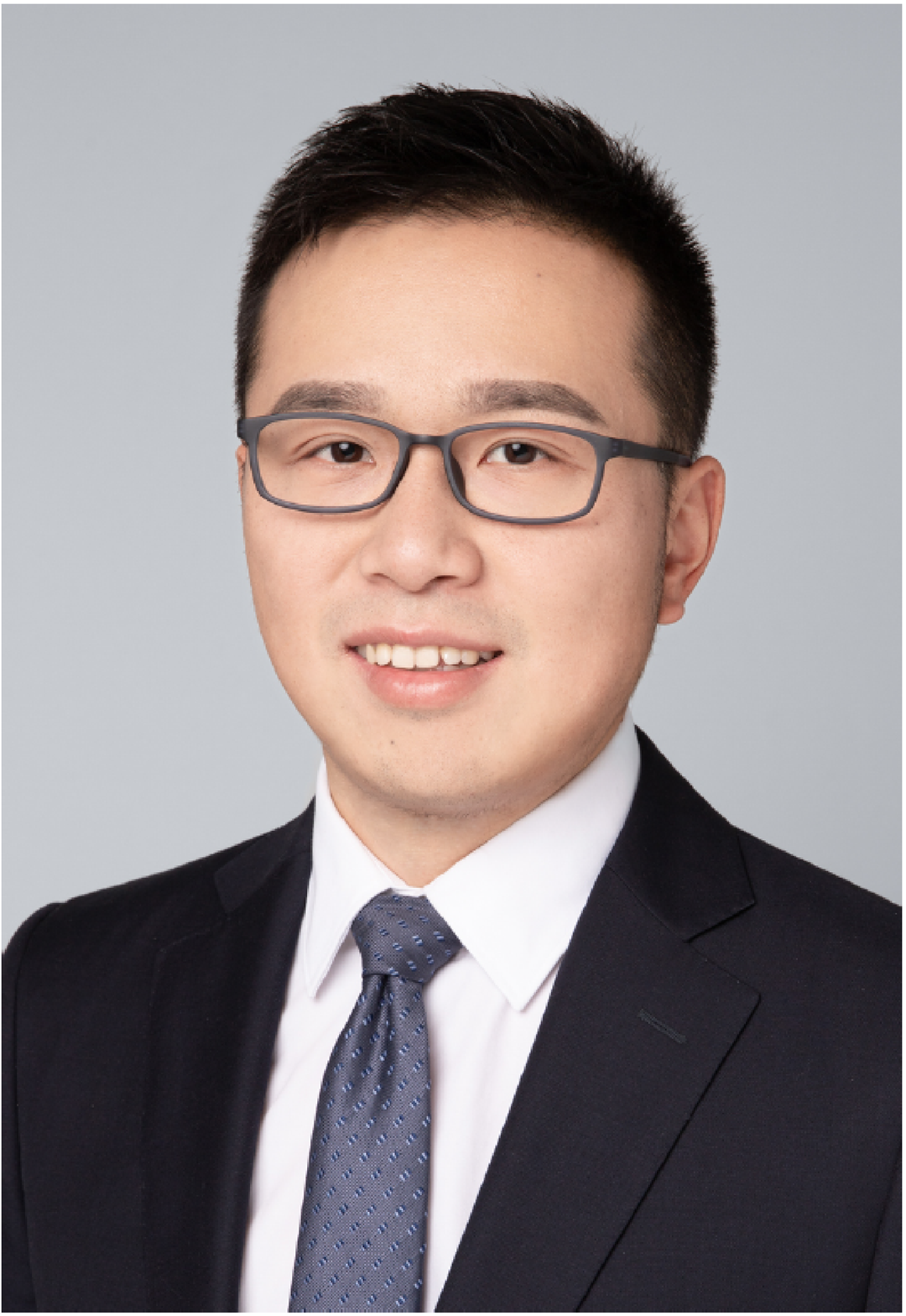}}]{Yang Tang} (Senior Member, IEEE) received the B.S. and Ph.D. degrees in electrical engineering from Donghua University, Shanghai, China, in 2006 and 2010, respectively. From 2008 to 2010, he was a Research Associate with The Hong Kong Polytechnic University, Hong Kong. From 2011 to 2015, he was a Post-Doctoral Researcher with the Humboldt University of Berlin, Berlin, Germany, and with the Potsdam Institute for Climate Impact Research, Potsdam, Germany. Since 2015, he has been a Professor with the East China University of Science and Technology, Shanghai. His current research interests include distributed estimation/control/optimization, cyber-physical systems, hybrid dynamical systems, computer vision, reinforcement learning and their applications.
	
	Prof. Tang was a recipient of the Alexander von Humboldt Fellowship and the ISI Highly Cited Researchers Award by Clarivate Analytics from 2017. He is a Senior Board Member of Scientific Reports, an Associate Editor of IEEE Transactions on Neural Networks and Learning Systems, IEEE Transactions on Emerging Topics in Computational Intelligence, IEEE Transactions on Circuits and Systems I: Regular Papers and IEEE Systems Journal, etc.
	
\end{IEEEbiography}









\end{document}